\documentclass[a4paper,fleqn]{cas-sc}

\usepackage[numbers]{natbib}

\usepackage{amssymb}
\usepackage{amsmath}
\usepackage{algorithm}
\usepackage{siunitx}
\usepackage{algpseudocode}

\def\tsc#1{\csdef{#1}{\textsc{\lowercase{#1}}\xspace}}
\tsc{WGM}
\tsc{QE}

\begin{document}
\let\WriteBookmarks\relax
\def\floatpagepagefraction{1}
\def\textpagefraction{.001}

\shorttitle{DPG method for mircowave-heated flows}

\shortauthors{Marquis et al.}

\title [mode = title]{A discontinuous Petrov–Galerkin finite-element framework for the simulation of microwave-heated flows.}



%

\author[1]{Oreste Marquis}[orcid=0000-0003-4771-3801,linkedin=orestemarquis]



\ead{oreste.marquis@etud.polymtl.ca}


\credit{Conceptualization, Formal Analysis, Investigation, Methodology, Software, Validation, Visualization, Writing - Original Draft}

\affiliation[1]{organization={Chemical engineering High-performance Analysis, Optimization and Simulation (CHAOS) laboratory, Department of Chemical Engineering, Polytechnique Montreal},
  addressline={2500 Chem. de Polytechnique},
  city={Montreal},
  postcode={H3T 1J4},
  state={Quebec},
  country={Canada}}

\author[2]{Matthias Maier}


\ead{maier@tamu.edu}

\ead[url]{https://matthiasmaier.org/index.html}

\credit{Conceptualization, Methodology, Writing - Review \& Editing}

\affiliation[2]{organization={Department of Mathematics, Texas A\&M University},
  addressline={400 Bizzell St.},
  city={College Station},
  postcode={TX 77840},
  state={Texas},
  country={United-States}}

\author[1]{Bruno Blais}[orcid=0000-0001-6053-6542,linkedin=bruno-blais-b39b407a]

\cormark[1]

\ead{bruno.blais@polymtl.ca}

\ead[url]{https://www.chaos-lab.ca/}

\credit{Conceptualization, Funding Acquisition, Methodology, Project Administration, Resources, Software, Supervision, Writing - Review \& Editing}

\cortext[1]{Corresponding author}



\begin{abstract}
  We present a high-order multiphysics solver for the simulation of microwave-heated flows. The solver couples a discontinuous Petrov–Galerkin (DPG) finite element method for the time-harmonic Maxwell equations with continuous Galerkin finite element methods for the heat equation and the incompressible Navier–Stokes equations. We validate the electromagnetic solver against multiple benchmark problems: wave propagation in a rectangular waveguide, a cavity problem with a singular solution, and a microwave-heated obstacle problem, comparing our results against numerical and experimental data from the literature. The results confirm the validity of the implementation and demonstrate its ability to perform adaptive mesh refinement using the DPG method's built-in error estimator. The final part of the study showcases the capabilities of the multiphysics framework through simulations of microwave-heated flow around obstacles with singular geometric features. These results highlight the potential of the proposed framework for the simulation and optimization of microwave-assisted chemical processes. Finally, the developed high-order multiphysics solver has a low memory footprint, since the electromagnetic solver relies on a Conjugate Gradient (CG) iterative solver and the fluid solver is implemented in a matrix-free fashion, making the overall approach scalable and well-suited for large-scale parallel simulations.
\end{abstract}


\begin{highlights}
  \item DPG implementation guide for time-harmonic Maxwell equations
  \item Low-memory and scalable multiphysics solver for microwave-heated flows
  \item Verification of the electromagnetic solver and its coupling using benchmark problems
\end{highlights}

\begin{keywords}
  Microwave Assisted Processes \sep Finite Element Method \sep Multiphysics Simulation \sep Heat Transfer \sep Time-Harmonic Maxwell Equations \sep Discontinuous Petrov–Galerkin
\end{keywords}

\maketitle

\section{Introduction}

Reducing the carbon footprint of the chemical industry is one of the central challenges of the energy transition \citep{IEA_NetZeroRoadmap_2023,IEA_TrackingCleanEnergy_2023}. The sector accounts for approximately 5\% of global CO$_2$ emissions, driven by its dependence on high-temperature operations and carbon-based feedstocks \citep{gabrielliNetzeroEmissionsChemical2023}. Since modern technologies increasingly depend on chemical products throughout their supply chains, simply cutting down chemical production is not a viable path forward; instead, the industry must pursue process innovation to reduce its environmental impact.\par

One promising avenue is microwave-assisted processes, first popularized in the chemical sciences by \citet{giguereApplicationCommercialMicrowave1986} and \citet{gedyeUseMicrowaveOvens1986}, in which direct molecular-level energy deposition enables selective and volumetric heating that can accelerate reaction kinetics, improve product selectivity, and reduce the formation of undesired by-products \citep{goyalReviewMicrowaveassistedProcess2022}. In addition, when powered by renewable energy sources, microwave-assisted processes offer a pathway toward the decarbonization of chemical processes. Despite broad applicability to drying, pyrolysis, sintering, and chemical synthesis \citep{valverdeStateArtFundamental2024}, industrial adoption remains limited by the sensitivity of electromagnetic fields to material properties and geometry, as well as by the inherent difficulty of measuring field quantities inside microwave reactors \citep{goyalReviewMicrowaveassistedProcess2022,yangContinuousFlowMicrowaveReactor2023}. These factors make the scale-up of microwave reactors a complex and costly endeavor.\par

To address these limitations, the scientific community has increasingly turned to numerical simulation tools, such as COMSOL Multiphysics \citep{COMSOL_Multiphysics_2026} and ANSYS HFSS \citep{ANSYS_HFSS_2026} based on the finite element method (FEM), to characterize and optimize microwave reactors \citep{goyalReviewMicrowaveassistedProcess2022,valverdeStateArtFundamental2024,tangReviewMultiphysicsSimulation2026}. The FEM is also what is adopted throughout the present work. \par

In all cases found in the chemical engineering literature, the electromagnetic field is solved in a time-harmonic formulation, which is justified by the fact that the microwave period can be shorter than millions of times the characteristic timescales of the heat transfer and fluid dynamics. This formulation is numerically challenging since accurately capturing the wave requires that both the mesh size $h$ and the polynomial degree $p$ be chosen appropriately to control the error $\tilde{e}$ \citep{ihlenburgDispersionAnalysisError1995a, melenkWavenumberExplicitConvergence2011a}:
\begin{equation}
  \tilde{e} \lesssim C_1
  \left[\left(\frac{hk}{2p}\right)^m + C_1 C_2 k\left(\frac{hk}{2p}\right)^{p + m} \right]
\end{equation}
where $k$ is the wavenumber, $m = \text{min}(l,p)$ with $l+1$ being the regularity of the solution, and $C_1$, $C_2$ are constants that are independent of $h$ and $k$ but depend on $p$. As a result, finer meshes and higher polynomial degrees are required as the wavenumber grows, rapidly inflating the number of degrees of freedom and the computational cost. It follows that multifrontal direct solvers become impractical at large scale due to their $\mathcal{O}(N^2)$ floating-point operations and $\mathcal{O}(N^{4/3})$ memory requirements \citep{petridesAdaptiveMultigridSolver2021a}, motivating the use of Krylov methods such as GMRES. However, effective preconditioning for those iterative solvers applied to the time-harmonic wave equation remains an open challenge: the linear system resulting from the discretization of the equation is indefinite, rendering classical approaches such as shifted Laplacian, multigrid, and domain decomposition methods less effective \citep{grahamNumericalAnalysisMultiscale2012a}. Resonance modes and material discontinuities further complexify the design of an efficient preconditioner and hinder convergence. \par

\subsection{Discontinuous Petrov-Galerkin method}
The challenges outlined above call for a discretization framework that is simultaneously stable, accurate under high-wavenumber conditions, and compatible with efficient iterative solvers. The discontinuous Petrov–Galerkin (DPG) method introduced by Demkowicz and Gopalakrishnan in a series of articles \citep{demkowiczClassDiscontinuousPetrov2010, demkowiczClassDiscontinuousPetrov2011, demkowiczClassDiscontinuousPetrov2012, zitelliClassDiscontinuousPetrov2011} addresses these challenges. The method yields several properties that make it particularly well-suited for time-harmonic problems. Most notably, it produces a Hermitian positive definite linear system, which enables the use of a Conjugate Gradient (CG) solver. Compared to alternatives such as GMRES, CG carries a significantly lower memory footprint as it does not require the storage in memory of the Krylov basis vectors. This property is critical for large-scale 3D simulations. Furthermore, the DPG method offers "absolute stability", meaning it does not require the mesh size $h$ to be "sufficiently small" for the convergence estimate to hold \citep{carstensenBreakingSpacesForms2016}.

Although wavenumber-explicit estimates of the method have only been mathematically established for the acoustic case \citep{demkowiczWavenumberExplicitAnalysis2012}, the structural analogy between the Helmholtz and Maxwell systems, combined with the numerical evidence presented for hexahedral and tetrahedral meshes \citep{carstensenBreakingSpacesForms2016}, gives confidence that similar behavior holds for Maxwell equations and that DPG remains a promising approach for the simulation of microwave heating processes. In addition, the stability of the DPG method allows for a robust error control mechanism using the built-in a posteriori error estimator, which can drive adaptive mesh refinement and maintain solution accuracy even in the presence of challenging features such as resonant modes, material interface discontinuities, or geometric singularities \citep{demkowiczDiscontinuousPetrovGalerkin2025}.\par

The primary software platform that has implemented DPG developments is the open-source package \texttt{hp3D} \citep{hp2024joss}, a Fortran-based framework supporting $hp$-adaptive discretization on mixed element meshes (tetrahedra, hexahedra, prisms, and pyramids). \texttt{hp3D} has demonstrated its versatility in some coupled multiphysics settings, including heat transfer and time-harmonic Maxwell equations, notably for the simulation of nonlinear thermal effects in optical fiber amplifiers \citep{hennekingVectorialEnvelopeMaxwell2025}.

\subsection{Novelty}

The DPG method is a mature finite element framework \citep{demkowiczDiscontinuousPetrovGalerkin2025}. Accordingly, the novelty of this work does not lie in the theoretical development of the DPG method itself, but in the following contributions:

\begin{itemize}
  \item \textbf{A practical implementation guide for time-harmonic wave problems.} While the mathematical description of the DPG method is well established in the literature, practical guides for its numerical implementation remain scarce, restricting its adoption to specialists in the field. This work addresses that gap by providing a detailed implementation guide for the electromagnetic wave equation, making the method more accessible to a broader community. It is complemented by a tutorial on the implementation of the DPG method for the acoustic wave equation\citep{marquisDealiiTutorialStep1002026a}, already added to the \texttt{deal.II} finite element library as part of the 9.8 release \citep{2026:arndt.bangerth.ea:deal}.

  \item \textbf{Integration of a DPG time-harmonic electromagnetic solver within a multiphysics framework for the simulation of microwave-heated flows.}  This implementation is carried out inside \texttt{lethe} \citep{alphoniusLethe10Opensource2026}, a validated multiphysics and multiphase solver targeted at chemical and manufacturing processes. \texttt{lethe} supports $h$-adaptivity and offers advanced matrix-free capabilities that enable large-scale parallel simulations \citep{prietosaavedraMatrixfreeStabilizedSolver2025}. The discretization order for each physics can be selected independently of one another, enabling a more flexible allocation of computational resources. This integration paves the way for coupled simulations of fluid dynamics, heat transfer, and time-harmonic electromagnetism. To the best of our knowledge, this work presents the first integration of the DPG method within an open-source multiphysics framework outside of the DPG community.
\end{itemize}

\subsection{Outline}

This work is organized as follows: Section~\ref{sec:dpg} provides a concise overview of the DPG method and its application to the time-harmonic Maxwell equations. Section~\ref{sec:multiphysics} describes the governing equations for heat transfer and fluid dynamics and discusses how they are coupled with the electromagnetic solver to form a fully coupled microwave-heated flow model. Section~\ref{sec:tests} presents numerical tests that verify the electromagnetic solver implementation and demonstrate the capabilities of the multiphysics coupling. Finally, Section~\ref{sec:conclusions} summarizes the key findings and discusses the current limitations of the framework alongside directions for future work.

\section{The DPG method for time-harmonic Maxwell equations}
\label{sec:dpg}

We present a brief introduction to the DPG method and its application to the time-harmonic Maxwell equations. For more details on the method, we refer the reader to the review of \citet{demkowiczDiscontinuousPetrovGalerkin2025}.

\subsection{Abstract problem formulation}

\begin{figure}
  \centering
  \includegraphics[width=0.75\textwidth]{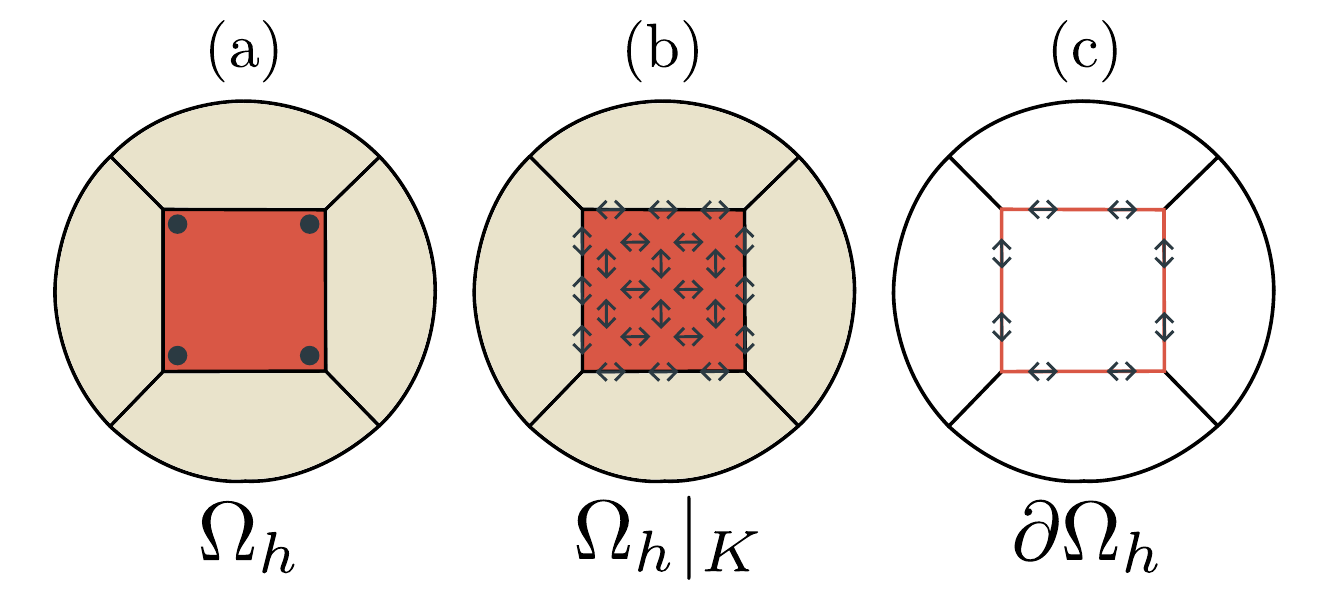}
  \caption{Schematic of the finite element spaces involved in the DPG formulation of the time-harmonic Maxwell equations on a two-dimensional circular domain $\Omega$ discretized using quadrilateral elements, together with its associated mesh skeleton $\partial\Omega_h$. The red cell and boundaries highlight a single element and its local degrees of freedom. (a) Discontinuous interior trial space $\mathcal{Q}^-_2\Lambda^2(\square_2)$, (b) discontinuous Nédélec test space $\mathcal{Q}^-_3\Lambda^1(\square_2)$ with polynomial degree one higher than the trial space, and (c) Nédélec trial trace space $\mathcal{Q}^-_2\Lambda^1(\square_2)$. Degrees of freedom are shown in dark grey, with arrows indicating vector-valued degrees of freedom. The notation $\Omega_h|_K$ emphasizes that the test space is discontinuous between each element $K$. } \label{fig:schematic_DPG}
\end{figure}

The DPG method belongs to the class of Petrov–Galerkin finite element approaches, in which the trial and test spaces are chosen independently of each other. The starting point is therefore an abstract variational problem posed on a trial space $U$ and a test space $V$, both assumed to be Hilbert spaces. Given a continuous bilinear form $b(\cdot,\cdot) : V \times U \to \mathbb{C}$, which stems from the weak formulation of the governing partial differential equations, and a continuous linear form $l(\cdot) : V \to \mathbb{C}$, which collects the source terms and the boundary data, the abstract problem reads:
\begin{equation}\label{eq:abstract_problem}
  \text{Find } u \in U \text{ such that } \quad b(v,u) = l(v), \quad \forall v \in V.
\end{equation}
The problem \eqref{eq:abstract_problem} is assumed to be well-posed, i.e., the bilinear form $b(\cdot,\cdot)$ satisfies an inf-sup condition on $V \times U$, so that the solution $u$ exists, is unique, and depends continuously on the data $l$. The DPG method discretizes \eqref{eq:abstract_problem} by relying on a mesh $\Omega_h$ consisting of elements $K$ that partition the domain $\Omega$, and whose element interfaces form the mesh skeleton $\partial\Omega_h = \bigcup_{K\in\Omega_h}\partial K$. A schematic of a discretization is presented in Figure \ref{fig:schematic_DPG}.

The practical formulation, which is the one implemented in this work, employs three finite-dimensional subspaces: a discontinuous (or broken) trial space $U_h$ corresponding to $U$, a discontinuous enriched\footnote{By enriched, it is meant that the finite element approximation order of the test space is of a higher degree than the one of the trial space.} test space $V_r$ corresponding to $V$, and a trace space $\hat U_h$ on the mesh skeleton that weakly enforces continuity between neighboring elements. Its unknowns are the discrete field $u_h$, the discrete trace $\hat u_h$, and the error representation function $\Psi \in V_r$, which represents the residual of the discrete solution in the test space. The formulation computes $\Psi$ while enforcing its orthogonality to $U_h$ and $\hat U_h$, the unknowns $u_h$ and $\hat u_h$ acting as the Lagrange multipliers of these constraints. It reads:
\begin{align}\label{eq:final_dpg_formulation}
   & \text{Find $u_h \in U_h$, $\hat{u}_h \in \hat{U}_h$ and $\Psi \in V_r$ such that}\nonumber             \\
   & \qquad \begin{cases}
              (v,\Psi)_{V_r} + b_h(v,u_h) + \langle v, \hat{u}_h \rangle_h = l_h(v), \quad \forall v \in V_r, \\
              b_h(w_h, \Psi) = 0, \quad \forall w_h \in U_h,                                                  \\
              \langle \hat{w}_h, \Psi \rangle_h = 0, \quad \forall \hat{w}_h \in \hat{U}_h.
            \end{cases}
\end{align}
Here, $b_h(\cdot,\cdot)$ and $l_h(\cdot)$ are the discrete counterparts of the bilinear form $b(\cdot,\cdot)$ and of the linear form $l(\cdot)$ of \eqref{eq:abstract_problem}, obtained by evaluating them element-wise on $\Omega_h$, which is why they depend on the mesh, $\langle \cdot,\cdot \rangle_h$ is the equally mesh-dependent duality pairing that couples the test functions to the interface unknowns across the skeleton $\partial\Omega_h$, and $(\cdot,\cdot)_{V_r}$ is the test-space inner product. All of the forms and pairings are defined element-wise on $\Omega_h$ and \textbf{are complex conjugate in the first argument}.

\subsection{Model problem and ultraweak formulation}

To simulate the electromagnetic field in microwave-heated processes, we begin with the dimensionless time-harmonic Maxwell equations\footnote{The convention used to make the system dimensionless is presented in Appendix \ref{appendix:dimensionless}.}:
\begin{subequations}\label{eq:maxwell_dimensionless}
  \begin{align}
     & \nabla \times \mathbf{E} - i \omega \mu_r \mathbf{H}                      = 0,          &  & (\text{ Faraday's law})        \\
     & \nabla \times \mathbf{H} + i \omega \varepsilon_{r,\text{eff}} \mathbf{E} = \mathbf{J}, &  & (\text{ Ampère-Maxwell's law})
  \end{align}
\end{subequations}
where $\mathbf{E}$ and $\mathbf{H}$ denote the electric and magnetic fields, $\omega$ is the angular frequency, $\varepsilon_{r,\text{eff}}$ is the effective relative permittivity defined as $\varepsilon_{r,\text{eff}} = \varepsilon_r + i \sigma_r$, which incorporates the reduced conductivity $\sigma_r$, $\mu_r$ is the relative permeability, and $\mathbf{J}$ is the current density. These equations are obtained using the time-harmonic ansatz
$\mathbf{\Phi}_\mathrm{total}(\mathbf{x},t) = \Re(\mathbf{\Phi}(\mathbf{x})e^{-i\omega t})$, where $\mathbf{\Phi}$ denotes an arbitrary field, $\mathbf{x}$ the position vector and $t$ the time. The equations \eqref{eq:maxwell_dimensionless} are supplemented with the following boundary conditions:
\begin{subequations}
  \begin{align}
     & \mathbf{n} \times \mathbf{E} = \mathbf{n} \times \mathbf{E}_D
     &                                                                                         & \text{on } \Gamma_D, \label{eq:dirichlet_boundary} \\
     & \mathbf{n} \times \mathbf{H} = \mathbf{n} \times \mathbf{H}_N - \mathbf{J}_{\mathrm{s}}
     &                                                                                         & \text{on } \Gamma_N, \label{eq:neumann_boundary}   \\
     & \mathbf{n} \times \mathbf{H}
    + Z_{\mathrm{s}}^{-1}\mathbf{n} \times (\mathbf{E}\times \mathbf{n})
    = \mathbf{g}
     &                                                                                         & \text{on } \Gamma_R. \label{eq:robin_boundary}
  \end{align}
\end{subequations}
In the above, $\Gamma_D$, $\Gamma_N$ and $\Gamma_R$ are the Dirichlet, Neumann and Robin boundaries, respectively, which form a disjoint partition $\Gamma = \Gamma_D \cup \Gamma_N \cup \Gamma_R$ of the domain boundary, $\mathbf{n}$ is the unit outward normal vector to the boundary, $Z_{\mathrm{s}}$ is the surface impedance of the Robin boundary, $\mathbf{E}_D$ and $\mathbf{H}_N$ are the prescribed electric and magnetic fields on the Dirichlet and Neumann boundaries, respectively, $\mathbf{J}_{\mathrm{s}}$ is a surface current density on the Neumann boundary, and $\mathbf{g}$ is a source term on the Robin boundary.\par

For the FEM discretization of the system of equations \eqref{eq:maxwell_dimensionless}, we use the ultraweak formulation rather than the more conventional primal form. We recall that the magnetic field is eliminated in the primal form to obtain a second-order equation for the electric field. The choice of the ultraweak formulation is motivated by the fact that it is useful to obtain both electric and magnetic field solutions to compute the electromagnetic power dissipation, as some materials can have a non-negligible magnetic loss. Additionally, the ultraweak formulation has the lowest inter-element regularity requirements and has shown good performance for high-frequency wave problems \citep{zitelliClassDiscontinuousPetrov2011,demkowiczWavenumberExplicitAnalysis2012,petridesAdaptiveMultigridSolver2021a,nagaraj3DDPGMaxwell2019,HENNEKING202130}. To obtain the ultraweak formulation, both equations presented in \eqref{eq:maxwell_dimensionless} are weakened by multiplying them with test functions $\mathbf{I}$ and $\mathbf{F}$, respectively, and integrating by parts. Using this ultraweak formulation of Maxwell's equations, and collecting the field unknowns in $u = (\mathbf{E},\mathbf{H})$, the trace unknowns in $\hat{u} = (\hat{\mathbf{E}},\hat{\mathbf{H}})$ and the test functions in $v = (\mathbf{F},\mathbf{I})$, we derive the mesh-dependent bilinear form $b_h$, the skeleton duality pairing $\langle \cdot,\cdot \rangle_h$ and the linear form $l_h$ that enter the DPG formulation \citep{hennekingNumericalStudyPollution2021,petridesAdaptiveMultigridSolver2021a}:
\begin{subequations}\label{eq:ultraweak_forms}
  \begin{align}
     & b_h(v,u)                      := (\nabla \times\mathbf{I},\mathbf{E})_{\Omega_h} - ( \mathbf{I}, i\omega \mu_r\mathbf{H})_{\Omega_h}
    + (\nabla \times \mathbf{F},\mathbf{H})_{\Omega_h} + (\mathbf{F}, i \omega \varepsilon_{r,\text{eff}} \mathbf{E})_{\Omega_h} - \langle \mathbf{F} ,Z_\mathrm{s}^{-1} \mathbf{n} \times (\mathbf{E} \times \mathbf{n})\rangle_{ \Gamma_R }, \label{eq:bilinear_form} \\
     & \langle v, \hat{u} \rangle_h  := \langle \mathbf{I}, \mathbf{n}\times \hat{\mathbf{E}} \rangle_{\partial \Omega_h} + \langle \mathbf{F} , \mathbf{n}\times \hat{\mathbf{H}}\rangle_{\partial \Omega_h \backslash  \Gamma_R }, \label{eq:trace_pairing}           \\
     & l_h(v)                        := (\mathbf{F},\mathbf{J})_{\Omega_h} -  \langle \mathbf{F} , \mathbf{g} \rangle_{ \Gamma_R }. \label{eq:linear_form}
  \end{align}
\end{subequations}
Here, $(\cdot,\cdot)_{\Omega_h}$ denotes the element-wise $L^2$ complex inner product, again with the complex conjugation in the first argument. The Robin boundary condition is applied to the electric field without loss of generality.

The forms \eqref{eq:ultraweak_forms} are defined on the following trial and test spaces, with $U := \mathcal{X}_{\mathrm{em}} \times \mathcal{X}_{\mathrm{em}}$, $\hat{U} := \hat{\mathcal{X}}_{\mathrm{em}} \times \hat{\mathcal{X}}_{\mathrm{em}}$ and $V := \mathcal{Y}_{\mathrm{em}} \times \mathcal{Y}_{\mathrm{em}}$:
\begin{align}
  \begin{cases}
     & \mathbf{E},\mathbf{H} \in \mathcal{X}_{\mathrm{em}}, \text{ with }\mathcal{X}_{\mathrm{em}}                          := (L^2(\Omega))^3,                                                      \\
     & \hat{\mathbf{E}},\hat{\mathbf{H}} \in \hat{\mathcal{X}}_{\mathrm{em}}, \text{ with }\hat{\mathcal{X}}_{\mathrm{em}}  := \bigl\{ \hat{\mathbf{x}} \in H^{-1/2}(\mathrm{curl},\partial\Omega_h)
    \;:\;
    \mathbf{n}\times \hat{\mathbf{x}} = \mathbf{n}\times \hat{\mathbf{x}}_\mathrm{D}
    \text{ on } \Gamma_D
    \bigr\},                                                                                                                                                                                         \\
     & \mathbf{F},\mathbf{I} \in \mathcal{Y}_{\mathrm{em}}, \:\:  \text{ with } \mathcal{Y}_{\mathrm{em}}                        :=  H(\mathrm{curl},\Omega_h),
  \end{cases}
\end{align}

where we refer the reader to Appendix \ref{appendix:polynomial_spaces} for the definition of the broken space $H(\mathrm{curl},\Omega_h)$, of the trace space $H^{-1/2}(\mathrm{curl},\partial\Omega_h)$, and more details on functional spaces or trace operators related to DPG in general. For an ultraweak formulation, note that the Neumann boundary condition is in fact a strongly imposed boundary condition on the flux variable (i.e., the magnetic field $\mathbf{H}$ in our current interpretation). This is why we have only stated the essential boundary condition for the trace space $\hat{\mathcal{X}}_{\mathrm{em}}$, which encompasses both the electric and magnetic fields, since we can directly impose the relations \eqref{eq:dirichlet_boundary} and \eqref{eq:neumann_boundary} on the relevant electric or magnetic field unknowns.

To complete the DPG discrete formulation of the problem, we define the inner product with which the test space is equipped, and the norm it induces, which determines the norm in which the residual is minimized. For the problem at hand, we use the \emph{adjoint graph norm} induced by the bilinear form $b_h(\cdot,\cdot)$ \citep{carstensenBreakingSpacesForms2016,petridesAdaptiveMultigridSolver2021a,hennekingNumericalStudyPollution2021}. For two test functions $v = (\mathbf{F},\mathbf{I})$ and $\tilde{v} = (\tilde{\mathbf{F}},\tilde{\mathbf{I}})$, the norm is defined as $\|v\|^2_{V_r} := (v,v)_{V_r}$, where
\begin{multline}\label{eq:test_inner_product}
  (v,\tilde{v})_{V_r} := \bigl( \nabla \times \mathbf{F} + i \omega \mu_r^* \mathbf{I},\;
  \nabla \times \tilde{\mathbf{F}} + i \omega \mu_r^* \tilde{\mathbf{I}} \bigr)_{\Omega_h}  + \bigl( \nabla \times \mathbf{I} - i \omega \varepsilon_{r,\mathrm{eff}}^* \mathbf{F},\;
  \nabla \times \tilde{\mathbf{I}} - i \omega \varepsilon_{r,\mathrm{eff}}^* \tilde{\mathbf{F}} \bigr)_{\Omega_h}               \\
  + \alpha \bigl( (\mathbf{F},\tilde{\mathbf{F}})_{\Omega_h} + (\mathbf{I},\tilde{\mathbf{I}})_{\Omega_h} \bigr)
  + \bigl( \mathbf{n} \times \mathbf{I} + (Z_\mathrm{s}^*)^{-1} \mathbf{n} \times (\mathbf{F} \times \mathbf{n}),\;
  \mathbf{n} \times \tilde{\mathbf{I}} + (Z_\mathrm{s}^*)^{-1} \mathbf{n} \times (\tilde{\mathbf{F}} \times \mathbf{n}) \bigr)_{\Gamma_R}.
\end{multline}
Here again, $(\cdot,\cdot)_{\Gamma_R}$ denotes the $L^2(\Gamma_R)$ complex inner product in its first argument, and the $^*$ symbole denotes complex conjugation. The constant $\alpha\in \mathcal{O}(1)$ is a strictly positive constant introduced to ensure that the test norm is localizable when the test space is broken. This is required for the element-wise computation of the error representation function in \eqref{eq:final_dpg_formulation} (we choose $\alpha=1$). The boundary contribution in the norm denoted by $\Gamma_R$ is used to enforce the Robin boundary condition and is required to obtain the convergence of the method when this type of boundary condition is present. Finally, the test norm induces the energy norm on the trial space as follows:
\begin{equation}\label{eq:energy_norm}
  \|u\|_E := \sup_{v \in V \setminus \{0\}} \frac{|b_h(v,u)|}{\|v\|_{V_r}},
\end{equation}
which is the norm in which the DPG solution is a best approximation \citep{demkowiczDiscontinuousPetrovGalerkin2025}, and in which the discretization error is measured in Section \ref{sec:tests}.

\subsection{Numerical implementation and discretization}

Our implementation, which is based on the \texttt{deal.II} library \citep{2026:arndt.bangerth.ea:deal}, discretizes the ultraweak formulation of the time-harmonic Maxwell equations using hexahedral finite elements. To exploit the broader range of features available in the \texttt{deal.II} library for real-valued problems, the trial spaces (i.e. $\mathbf{E}, \mathbf{H}, \hat{\mathbf{E}}, \hat{\mathbf{H}}$) and test functions (i.e. $\mathbf{F}, \mathbf{I}$) are decomposed into their real (e.g. $\mathbf{E}_{\mathrm{re}}$) and imaginary (e.g. $\mathbf{E}_{\mathrm{im}}$) parts, and the resulting system is solved entirely using regular double precision floating points data types.

The discrete trial, trace, and test spaces are defined as follows:
\begin{align}
  \begin{cases}
     & U_h = \mathbf{E}_\mathrm{re} \times \mathbf{E}_\mathrm{im} \times \mathbf{H}_\mathrm{re} \times \mathbf{H}_\mathrm{im},                              \\
     & \hat{U}_h = \hat{\mathbf{E}}_\mathrm{re} \times \hat{\mathbf{E}}_\mathrm{im} \times \hat{\mathbf{H}}_\mathrm{re} \times \hat{\mathbf{H}}_\mathrm{im} \\
     & V_r = \mathbf{F}_\mathrm{re} \times \mathbf{F}_\mathrm{im} \times \mathbf{I}_\mathrm{re} \times \mathbf{I}_\mathrm{im}.
  \end{cases}
\end{align}
Each component of the interior trial space $u_h \in U_h$ is discretized using discontinuous finite elements $\mathcal{Q}^-_p \Lambda^3(\square_3)$ constructed from tensor products of Lagrange polynomials of degree $p$ in each spatial direction, yielding a discrete approximation of $L^2(\Omega_h)$. The trace unknowns $\hat{u}_h \in \hat{U}_h$ are discretized using Nédélec elements of the first kind $\mathcal{Q}^-_p \Lambda^1(\square_3)$ of degree $p$. The tangential trace operator $\mathrm{tr}_{\mathrm{curl},\top}$ is applied to these elements, and the interior degrees of freedom are discarded to ensure that the resulting finite element space provides a conforming discretization of $H^{-1/2}(\mathrm{curl}, \partial\Omega_h)$ functions. The test space $v \in V_r$ is also discretized using Nédélec elements of the first kind $\mathcal{Q}^-_p \Lambda^1(\square_3)$, but without applying the trace operator, so that the resulting space is conforming to $H(\mathrm{curl}, \Omega_h)$. In addition, the degree of the polynomials used for the test space is set to $p + \Delta p$, where $\Delta p \in \mathbb{N}^+$ is the enrichment parameter required by the method. In what follows, we use $\Delta p = 1$. Accordingly, the discrete trial and test functions are expanded in terms of the corresponding basis functions as:
\begin{align}\begin{cases}
     & \mathbf{u}_h^{(k)}(\mathbf{x}) = \sum_{i} w^{(k)}_{i} \boldsymbol{\phi}_{i}(\mathbf{x}), \\  & \hat{\mathbf{u}}_h^{(k)}(\mathbf{x}) = \sum_{i} \hat{w}^{(k)}_{i} \hat{\boldsymbol{\phi}}_{i}(\mathbf{x}), \\ & \mathbf{v}^{(k)}(\mathbf{x}) = \sum_{i} q^{(k)}_{i} \boldsymbol{\psi}_{i}(\mathbf{x}),
  \end{cases}
\end{align}
where $\boldsymbol{\phi}_i$, $\hat{\boldsymbol{\phi}}_i$ and $\boldsymbol{\psi}_i$ are the basis functions for the interior trial, trace and test spaces, respectively, and $w_i$, $\hat{w}_i$ and $q_i$ are the corresponding coefficients. The superscript $k\in \{1,2,3,4\}$ indexes the fields components of the element, for example $(\mathbf{E}_\mathrm{re}, \mathbf{E}_\mathrm{im}, \mathbf{H}_\mathrm{re}, \mathbf{H}_\mathrm{im})$. \par

The above discretization, supplemented by the abstract formulation \eqref{eq:final_dpg_formulation}, leads to the following linear system in a DPG framework, where each matrix entry implicitly carries a composite index $(i,k)$, where $i$ is the index of the finite element basis function and $k$ is the index of the field component:
\begin{equation}
  \begin{bmatrix}
    G               & B & \hat{B} \\
    B^\dagger       & 0 & 0       \\
    \hat{B}^\dagger & 0 & 0
  \end{bmatrix}
  \begin{bmatrix}
    \Psi \\
    u_h  \\
    \hat{u}_h
  \end{bmatrix}
  =
  \begin{bmatrix}
    l \\
    0 \\
    0
  \end{bmatrix},
\end{equation}
where the entries of the interior matrix $B$, of the interface matrix $\hat{B}$, of the Gram matrix $G$ and of the load vector $l$ are given by
\begin{align*}
  \begin{cases}
    B_{(i,k),(j,l)} = b_h( \boldsymbol{\psi}^{(k)}_i, \boldsymbol{\phi}^{(l)}_j),                         & \text{see } \eqref{eq:bilinear_form},      \\
    \hat{B}_{(i,k),(j,l)} = \langle \boldsymbol{\psi}^{(k)}_i, \hat{\boldsymbol{\phi}}^{(l)}_j \rangle_h, & \text{see } \eqref{eq:trace_pairing},      \\
    G_{(i,k),(j,l)} = (\boldsymbol{\psi}^{(k)}_i, \boldsymbol{\psi}^{(l)}_j)_{V_r},                       & \text{see } \eqref{eq:test_inner_product}, \\
    l_{(i,k)} = l_h(\boldsymbol{\psi}^{(k)}_i),                                                           & \text{see } \eqref{eq:linear_form}.
  \end{cases}
\end{align*}
In our framework, this system is condensed by eliminating the error representation function $\Psi$, resulting in the following:
\begin{equation}
  \begin{bmatrix}
    B^\dagger G^{-1}B       & B^\dagger G^{-1}\hat{B}       \\
    \hat{B}^\dagger G^{-1}B & \hat{B}^\dagger G^{-1}\hat{B}
  \end{bmatrix}
  \begin{bmatrix}
    u_h \\
    \hat{u}_h
  \end{bmatrix}
  =
  \begin{bmatrix}
    B^\dagger G^{-1}l \\
    \hat{B}^\dagger G^{-1}l
  \end{bmatrix}.
\end{equation}
Recall that this is achievable since the Gram matrix $G$ is block-diagonal and local to each element, due to the broken nature of the test space. Consequently, we follow the approach taken by \citep{petridesAdaptiveDPGMethod2017a}, and we perform another static condensation step to eliminate the interior unknowns $u_h$ and obtain a system only in terms of the trace unknowns $\hat{u}_h$. Defining the matrix blocks as $M_1 = B^\dagger G^{-1}B$, $M_2 = B^\dagger G^{-1}\hat{B}$, $M_3 = \hat{B}^\dagger G^{-1}\hat{B}$, $M_4 = B^\dagger G^{-1}$ and $M_5 = \hat{B}^\dagger G^{-1}$, the resulting system is:
\begin{equation}\label{eq:linear_system}
  (M_3 - M_2^\dagger M_1^{-1} M_2) \hat{u}_h = (M_5 - M_2^\dagger M_1^{-1} M_4) l.
\end{equation}
Equation \eqref{eq:linear_system} is the linear system that is solved for. It corresponds to the Schur complement associated with the elimination of the interior unknowns. After solving for the trace unknowns $\hat{u}_h$, the interior solution $u_h$ and the error representation function $\Psi$ are recovered locally on each element by solving:
\begin{equation}
  u_h = M_1^{-1} (M_4 l - M_2 \hat{u}_h), \quad \Psi = G^{-1} (B u_h + \hat{B} \hat{u}_h - l).
\end{equation}
The residual norm on each element is then computed as $\|\Psi\|_{V_r} = \sqrt{\Psi^\dagger G \Psi}$ and used as the error estimator in the adaptive mesh refinement strategy. A pseudocode for the electromagnetic solver is presented in Algorithm \ref{alg:solve_maxwell}, and one for the system assembly is presented in Algorithm \ref{alg:dpg_assembly}. \par

\begin{algorithm}[ht]
  \caption{Pseudocode for solving the time-harmonic Maxwell equations.}
  \label{alg:solve_maxwell}
  \begin{algorithmic}[1]
    \Function{SolveTimeHarmonicMaxwell}{ }
    \State $\Omega_h$ $\gets$ MakeDiscretizedDomain()
    \State $U_h,\hat{U}_h,V_r$ $\gets$ InitializeFiniteElementSpaces()
    \State $A,b$ $\gets$ AssembleLinearSystem($\Omega_h, U_h,\hat{U}_h, V_r$)
    \State $\hat{\mathbf{E}}_\mathrm{re},\hat{\mathbf{E}}_\mathrm{im}, \hat{\mathbf{H}}_\mathrm{re}, \hat{\mathbf{H}}_\mathrm{im}$ $\gets$ SolveConjugateGradient($A,b$)
    \State $\mathbf{E}_\mathrm{re}, \mathbf{E}_\mathrm{im}, \mathbf{H}_\mathrm{re},\mathbf{H}_\mathrm{im}$ $\gets$ ReconstructSolution($\hat{\mathbf{E}}_\mathrm{re},\hat{\mathbf{E}}_\mathrm{im}, \hat{\mathbf{H}}_\mathrm{re}, \hat{\mathbf{H}}_\mathrm{im}$)
    \If{MeshAdaptation()}
    \State $\Psi$ $\gets$ ReconstructResidual($\mathbf{E}_\mathrm{re}, \mathbf{E}_\mathrm{im}, \mathbf{H}_\mathrm{re},\mathbf{H}_\mathrm{im},\hat{\mathbf{E}}_\mathrm{re},\hat{\mathbf{E}}_\mathrm{im}, \hat{\mathbf{H}}_\mathrm{re}, \hat{\mathbf{H}}_\mathrm{im}$)
    \EndIf
    \EndFunction
  \end{algorithmic}
\end{algorithm}

\begin{algorithm}[ht]
  \caption{Pseudocode for the assembly of the DPG linear system.}
  \label{alg:dpg_assembly}
  \begin{algorithmic}[1]
    \Function{AssembleLinearSystem}{$\Omega_h, U_h,\hat{U}_h, V_r$}
    \ForAll{Cells $K \in \Omega_h$}
    \ForAll{Quadrature points $q \in K$}
    \State $G,B,l$ $\gets$ ComputeInteriorContribution($U_h,V_r$)
    \EndFor
    \ForAll{Faces $\partial K \in K$}
    \ForAll{Quadrature points $q \in \partial K$}
    \If{$q \in \Gamma_R$}
    \State $G,\hat{B},l$ $\gets$ ComputeRobinBoundaryContribution($\hat{U}_h, V_r$)
    \Else
    \State $\hat{B}$ $\gets$ ComputeSkeletonContribution($\hat{U}_h, V_r$)
    \EndIf
    \EndFor
    \EndFor
    \State $M_1, M_2, M_3, M_4, M_5$ $\gets$ BuildCondensationMatrix($G,B,\hat{B}$)
    \State $A_K$ $\gets$ BuildCellMatrix($M_1,M_2,M_3$)
    \State $b_K$ $\gets$ BuildCellRHS($M_1,M_2,M_4,M_5,l$)
    \State $A,b$ $\gets$ DistributeLocalToGlobal($A_K,b_K$)
    \EndFor
    \EndFunction
  \end{algorithmic}
\end{algorithm}

Finally, to solve the condensed linear system \eqref{eq:linear_system}, we employ a CG solver without preconditioning because it is sufficient for our current purpose. Additionally, preliminary numerical experiments indicate that standard black-box preconditioners, such as algebraic multigrid, SSOR, or incomplete LU factorization, do not lead in satisfactory improvement of the wall clock time. It follows that the implementation of a preconditioner tailored to the present formulation is left for future work. We also note that the current formulation of the Robin boundary condition, combined with the ultraweak formulation, leaves certain magnetic field degrees of freedom unconstrained at those boundaries. We have found that this does not affect convergence for $p \geq 1$ when using a Conjugate Gradient linear solver.

\section{Multiphysics coupling}
\label{sec:multiphysics}

To simulate microwave-heated flows, the electromagnetic DPG solver is implemented within \texttt{lethe} where it can be coupled with the heat transfer and fluid dynamics solvers already available. In this section, we describe briefly the implementation details of the additional physics equations we consider within \texttt{lethe} and how they couple with the electromagnetic solver to simulate microwave-heated flows. For more details on the implementation of the fluid dynamics and heat transfer solvers in \texttt{lethe}, we refer the reader to \citet{alphoniusLethe10Opensource2026}.

\subsection{Fluid dynamics}

At its core, \texttt{lethe} is a computational fluid dynamics software framework that solves the incompressible Newtonian Navier-Stokes equations:
\begin{subequations}
  \begin{align}
     & \frac{\partial \mathbf{u}}{\partial t} + (\mathbf{u} \cdot \nabla) \mathbf{u} = -\nabla p_{\mathrm{r}}  + \nu \nabla^2  \mathbf{u} \label{eq::single_phase_ns_momentum} \\
     & \nabla \cdot \mathbf{u} = 0    \label{eq::single_phase_ns_continuity}
  \end{align}
\end{subequations}
where $\mathbf{u}$ is the velocity, $\nu$ is the kinematic viscosity and $p_{\mathrm{r}}$ is the kinematic pressure ($p_{\mathrm{r}}=\frac{p}{\rho}$, with $p$ the pressure and $\rho$ the density). Using the scalar test function $q$ for the continuity equation and the vector test function $\mathbf{v}$ for the momentum equation, the above system is discretized, and the resulting weak form of the equations is:
\begin{subequations} \label{eq:weak_NS}
  \begin{align}
     & (\mathbf{v}, \frac{\partial \mathbf{u}}{\partial t})_{\Omega_h}  + (\mathbf{v}, (\mathbf{u} \cdot \nabla) \mathbf{u})_{\Omega_h} = (\nabla \cdot \mathbf{v}, p_{\mathrm{r}})_{\Omega_h} - \nu (\nabla \mathbf{v}, \nabla \mathbf{u})_{\Omega_h} \\
     & (q, \nabla \cdot \mathbf{u})_{\Omega_h}  = 0
  \end{align}
\end{subequations}
In \eqref{eq:weak_NS}, the products between spaces are the standard $L^2$ inner product generalized to tensors of any rank. The associated trial and test spaces are defined as:
\begin{align}
  \begin{cases}
     & \mathbf{u} \in \mathcal{X}_\mathbf{u}, \text{ with } \mathcal{X}_\mathbf{u} := \{\mathbf{u} \in (H^1(\Omega_h))^3 : \mathbf{u} = \mathbf{u}_D \text{ on } \Gamma_D\}, \\
     & p_{\mathrm{r}} \in \mathcal{X}_p, \text{ with } \mathcal{X}_p := L^2(\Omega_h),                                                                                       \\
     & \mathbf{v} \in \mathcal{Y}_\mathbf{v}, \text{ with } \mathcal{Y}_\mathbf{v} := \{\mathbf{v} \in (H^1(\Omega_h))^3 : \mathbf{v} = 0 \text{ on } \Gamma_D\},            \\
     & q \in \mathcal{Y}_q, \text{ with } \mathcal{Y}_q := L^2(\Omega_h).
  \end{cases}
\end{align}

To control numerical instabilities associated with the convection term in the momentum equation and the saddle-point nature of the system, the above equations are modified using the Streamline Upwind Petrov-Galerkin (SUPG) stabilization for the momentum equation and the Pressure Stabilizing Petrov-Galerkin (PSPG) stabilization for the continuity equation. Both stabilizations are implicit and based on a parameter that is computed based on the local flow conditions and mesh size. For further details on how the nonlinear system arising from Eq. \eqref{eq:weak_NS} is solved, or on the matrix-free implementation, we refer the reader to \citet{prietosaavedraMatrixfreeStabilizedSolver2025}.

\subsection{Heat transfer}

To simulate the heat transfer, \texttt{lethe} solves for enthalpy ($\mathcal{H}= c_\mathrm{p} T$) conservation assuming an incompressible flow:
\begin{equation} \label{eq:heat_transfer}
  \rho c_\mathrm{p}\frac{\partial T}{\partial t}+ \mathbf{u} \cdot \nabla \left( \rho c_\mathrm{p} T\right) =\nabla \cdot (\kappa \nabla T) + Q_\text{em} ,
\end{equation}
where $c_\mathrm{p}$ is the isobaric specific heat capacity, $T$ is the temperature, $\kappa$ is the thermal conductivity, and $Q_\text{em}$ represents the electromagnetic power dissipation. For time-harmonic electromagnetic fields, this is obtained from the time-averaged Poynting vector $\overline{\mathbf{S}}$:
\begin{equation}
  Q_\text{em}=-\nabla \cdot \overline{\mathbf{S}} = \frac{1}{2}\sigma|\mathbf{E}|^2 + \frac{1}{2}\omega\varepsilon_0\varepsilon_\mathrm{im}|\mathbf{E}|^2 + \frac{1}{2}\omega\mu_0\mu_\mathrm{im}|\mathbf{H}|^2,
\end{equation}
where $|\cdot|^2$ is the squared modulus of the complex field amplitudes, the overbar symbole denotes the time average, the parameters $\varepsilon_\mathrm{im}$ and $\mu_\mathrm{im}$ are the imaginary parts of the relative permittivity and permeability, respectively, and $\sigma$ is the conductivity. Although we consider temperature-independent physical properties in the present work, \texttt{lethe} is capable of simulating flows with temperature-dependent physical properties for the time-harmonic Maxwell equations, the fluid dynamics, and the heat transfer.

Since the DPG electromagnetic solver is implemented in a dimensionless form, the resulting electromagnetic field solutions need to be dimensionalized back before the electromagnetic power dissipation is computed. This is done either directly from user-input electric or magnetic field reference values (i.e., $E_0$ or $H_0 = E_0/Z_0 $, respectively) or from the input power at the waveguide port. The latter requires computing the dimensionless time-averaged Poynting vector at the waveguide port surface $A$ from the electromagnetic excitation and then using the result to match the user input power, i.e., by isolating the reference electric field $E_0$ in the following equation:
\begin{equation}
  \overline{P}_\text{input} = \frac{1}{2} \frac{E_0^2}{Z_0} \int_A \Re{(\mathbf{E} \times \mathbf{H}^*)} \cdot \mathbf{n} \text{d}A. \label{eq:power_scaling}
\end{equation}
The weak form of the enthalpy conservation equation \eqref{eq:heat_transfer} reads:
\begin{align}
  \  & (\theta, \rho c_\mathrm{p}\frac{\partial T}{\partial t})_{\Omega_h}  + (\theta, \mathbf{u}\cdot\nabla(\rho c_\mathrm{p}T))_{\Omega_h} =  -(\nabla \theta, \kappa\nabla T)_{\Omega_h} + (\theta, Q_\text{em})_{\Omega_h}, \label{eq:heat_transfer_weak_form}
\end{align}
with $\theta$ the test function. The associated trial and test spaces are defined as:
\begin{align}
  \begin{cases}
     & T \in \mathcal{X}_T, \text{ with } \mathcal{X}_T := \{T \in H^1(\Omega_h) : T = T_D \text{ on } \Gamma_D\},                        \\
     & \theta \in \mathcal{Y}_\theta, \text{ with } \mathcal{Y}_\theta := \{\theta \in H^1(\Omega_h) : \theta = 0 \text{ on } \Gamma_D\}.
  \end{cases}
\end{align}
Similarly to the fluid dynamics equations, Eq. \eqref{eq:heat_transfer_weak_form}
also uses SUPG stabilization to mitigate the numerical instabilities associated with high Péclet number flows.

\subsection{Time integration and coupling strategy}

The time harmonic assumption for the electromagnetic solver decouples the solution of the electromagnetic fields from the time integration of the fluid dynamics and heat transfer equations. This decoupling is physically justified by the large-scale separation between the characteristic electromagnetic response times and the fluid or thermal characteristic times. Indeed, in the microwave frequency range, the electromagnetic fields oscillate at periods of nanoseconds, while the characteristic time scales for fluid flow and heat transfer of a typical microwave-heated process are of the order of milliseconds for the fluid and seconds for the heating:
\begin{equation}
  \tau_{\mathrm{em}} \ll \tau_{\mathrm{adv}} \sim \frac{L_c}{U_c} \ll \tau_{\mathrm{heating}} \sim \frac{\rho c_\mathrm{p} \Delta T}{Q_\text{em}},
\end{equation}
where $L_c$ and $U_c$ are the characteristic length and velocity of the flow, respectively, and $\Delta T$ is the characteristic temperature change due to the microwave heating. Therefore, the electromagnetic fields reach their steady-state spatial distribution essentially instantaneously relative to any appreciable change in the material properties. Under this assumption, the time-averaged electromagnetic power dissipation $Q_\text{em}$ acts as a quasi-static volumetric heat source that evolves only through its dependence on the material properties, which themselves change on the slow heat transfer time scale.\par

This separation of scales motivates the following sequential coupling strategy: the electromagnetic fields are solved first at $t=0$ to obtain the initial power dissipation field, and then solved again if necessary (e.g., mesh adaptation, a change in $\varepsilon_{r,\mathrm{eff}}$ or $\mu_r$). Between electromagnetic solves, the power dissipation field $Q_\text{em}$ is frozen and reused as a fixed source term in the heat equation. In practice, this means the electromagnetic solver is called orders of magnitude less frequently than the fluid dynamics and heat transfer solvers, resulting in substantial computational savings for the overall coupled simulation. The pseudocode for this coupling strategy is summarized in Algorithm~\ref{alg:lethe_time_integration_em}.

\begin{algorithm}
  \caption{Pseudocode for the multiphysics coupling.}
  \label{alg:lethe_time_integration_em}
  \begin{algorithmic}[1]
    \Function{Integrate}{$\Delta t$, $t_{\mathrm{end}}$}
    \If{t == 0}
    \State $\mathbf{E}^{(0)}, \mathbf{H}^{(0)} \gets$ SolveTimeHarmonicMaxwell()
    \State $Q_{\mathrm{em}}^{(0)} \gets$ ComputePowerDissipation($\mathbf{E}^{(0)}, \mathbf{H}^{(0)}$)
    \State $\mathbf{u}^{(0)}, p_{\mathrm{r}}^{(0)} \gets$ SolveNavierStokes()
    \State $T^{(0)} \gets$ SolveHeatEquation($Q_{\mathrm{em}}^{(0)}, \mathbf{u}^{(0)}$)
    \State $\varepsilon_{r,\mathrm{eff}}^{(\Delta t)}, \mu_r^{(\Delta t)} \gets$ UpdateMaterialProperties($T^{(0)}$)
    \Else
    \While{ $t < t_{\mathrm{end}}$}
    \State $t \gets t + \Delta t$
    \If{ShouldUpdateEMFields($\varepsilon_{r,\mathrm{eff}}^{(t)}, \mu_r^{(t)}$)}
    \State $\mathbf{E}^{(t)}, \mathbf{H}^{(t)} \gets$ SolveTimeHarmonicMaxwell()
    \State $Q_{\mathrm{em}}^{(t)} \gets$ ComputePowerDissipation($\mathbf{E}^{(t)}, \mathbf{H}^{(t)}$)
    \Else
    \State $\mathbf{E}^{(t)}, \mathbf{H}^{(t)}, Q_{\mathrm{em}}^{(t)} \gets \mathbf{E}^{(t-\Delta t)}, \mathbf{H}^{(t-\Delta t)}, Q_{\mathrm{em}}^{(t-\Delta t )} $
    \EndIf
    \State $\mathbf{u}^{(t)}, p_{\mathrm{r}}^{(t)} \gets$ SolveNavierStokes()
    \State $T^{(t)} \gets$ SolveHeatEquation($Q_{\mathrm{em}}^{(t)}, \mathbf{u}^{(t)}$)
    \State $\varepsilon_{r,\mathrm{eff}}^{(t+\Delta t)}, \mu_r^{(t+\Delta t)} \gets$ UpdateMaterialProperties($T^{(t)}$)
    \EndWhile
    \EndIf
    \EndFunction
  \end{algorithmic}
\end{algorithm}

\section{Numerical tests}
\label{sec:tests}
In this section, we consider three cases to verify the implementation of the DPG electromagnetic solver. The first two cases target the electromagnetic solver alone, while the last case demonstrates the capabilities of the multiphysics coupling for microwave-heated flows. All of the cases are simulated using 3D hexahedral meshes, but with varying polynomial degree $p$ for the electromagnetic, fluid dynamics, and heat transfer problems. The simulations parameters for each case are summarized in Appendix \ref{appendix:simulation_parameters}. \par

\subsection{Waveguide analytical verification} \label{sec:waveguide_validation}

We consider the classical problem of a rectangular waveguide with perfectly conducting walls that propagates an electromagnetic excitation in a TE$_{10}$ mode traveling in the $x_3$-direction. This setup has a known analytical solution that can be used to compute the error and convergence rates of the method and verify the implementation. The dimensionless analytical solution of a TE$_{mn}$ mode that we use for this test case is given by:
\begin{align}
  \mathbf{E}_{\mathrm{sol}} & = i\frac{\omega \mu_r}{k_c^2}  \begin{bmatrix} -k_{x_2} \cos(k_{x_1} x_1) \sin(k_{x_2} x_2) \\ k_{x_1} \sin(k_{x_1} x_1) \cos(k_{x_2} x_2) \\ 0 \end{bmatrix} e^{i k_{x_3} x_3},                                                                                                                                                    \\
  \mathbf{H}_{\mathrm{sol}} & = \begin{bmatrix} - i \frac{k_{x_3} k_{x_1}}{k_c^2} \sin(k_{x_1} x_1) \cos(k_{x_2} x_2)   \\-i \frac{k_{x_3} k_{x_2}}{k_c^2} \cos(k_{x_1} x_1) \sin(k_{x_2} x_2)\\ \cos(k_{x_1} x_1) \cos(k_{x_2} x_2) \end{bmatrix} e^{i k_{x_3} x_3},
\end{align}
with $k_{x_1} = \frac{m \pi}{a}$, $k_{x_2} = \frac{n \pi}{b}$, $k_c^2 = k_{x_1}^2 + k_{x_2}^2$ and $k_{x_3} = \sqrt{\omega^2 \varepsilon_{r,\mathrm{eff}} \mu_r - k_c^2}$. The parameters $a$ and $b$ are the width and height of the waveguide, respectively.\par

\begin{figure}
  \centering
  \includegraphics[width=0.45\textwidth]{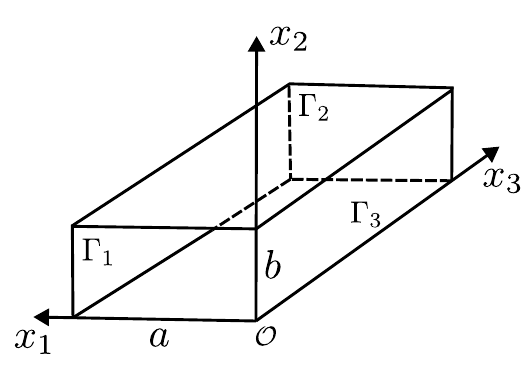}
  \caption{Schematic of the waveguide geometry.}\label{fig:waveguide_schematic}
\end{figure}

The parameters $a$, $b$, and the waveguide length are chosen so that the simulation domain matches the scale of a typical laboratory microwave reactor. Additionally, the excitation frequency is set to \SI{2.45}{\giga\hertz}, the standard frequency for microwave-heating applications. This way, by using the void properties for the material filling the waveguide (i.e., $\varepsilon_{r,\mathrm{eff}} =1$ and $\mu_r = 1$), we can ensure that the method is at least sufficient to capture the minimal dimensionless wavenumber $\tilde{\kappa}$ of relevant microwave heating problems. This number is defined as:
\begin{equation}
  \tilde{\kappa} = \omega \sqrt{\varepsilon \mu} L_c,
\end{equation}
and measures how many cycles of electromagnetic oscillations there are per unit of characteristic length. It is a key parameter that describes the frequency regime of the problem and dictates how challenging Maxwell's equations are numerically. Verifying that the method can at least accurately capture the solution for the minimal $\tilde{\kappa}$ value is therefore an essential baseline. With this in mind, the dimensions we choose are a width $a=\SI{0.25}{\metre}$, a height $b=\SI{0.25}{\metre}$ and a length $L=\SI{1}{\metre}$, which results in a dimensionless wavenumber $\tilde{\kappa} \approx 16 \pi$. A schematic of the geometry is shown in Figure \ref{fig:waveguide_schematic}. \par

In addition to the geometry and frequency, we also choose the boundary conditions for this test case to be representative of a typical waveguide setup, mirroring how electromagnetic fields are excited in our target applications. In particular, we use a port boundary condition at the inlet of the waveguide to excite the TE$_{10}$ mode, perfect electric conductor (PEC) boundary conditions on the walls, and an impedance matching boundary condition at the outlet to minimize the reflections. Those are defined as follows:
\begin{align}
   & \mathbf{n} \times \mathbf{H} + \frac{k_{x_3}}{\omega \mu_r} \mathbf{n} \times ( \mathbf{E} \times \mathbf{n} ) = \mathbf{n} \times \mathbf{H}_{\mathrm{TE}_{10}} + \frac{k_{x_3}}{\omega \mu_r} \mathbf{n} \times ( \mathbf{E}_{\mathrm{TE}_{10}} \times \mathbf{n} ) \quad \text{on } \Gamma_1 \\
   & \mathbf{n} \times \mathbf{H} + \frac{k_{x_3}}{\omega \mu_r} \mathbf{n} \times ( \mathbf{E} \times \mathbf{n} ) = 0 \quad \text{on } \Gamma_2,
  \\
   & \mathbf{n} \times \mathbf{E} = 0, \quad \text{on } \Gamma_3,
\end{align}
where $\Gamma_1$, $\Gamma_2$ and $\Gamma_3$ correspond to the inlet, outlet and walls of the waveguide respectively. \par

\begin{figure}
  \centering
  \includegraphics[width=0.5\textwidth]{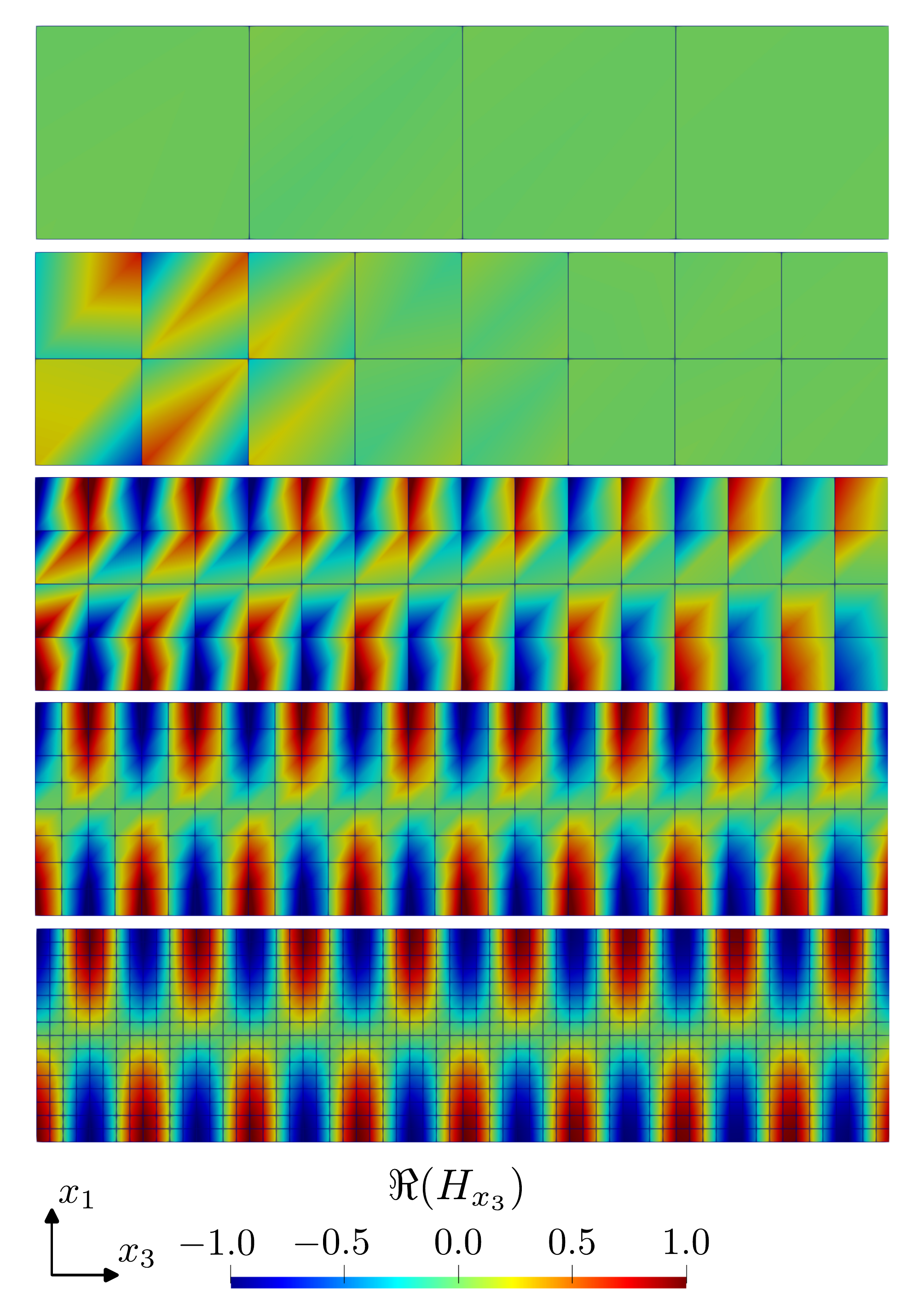}
  \caption{Real part of the $x_3$-component of the magnetic field for the waveguide test case at successive mesh refinement levels with $p=1$. The top panel shows the solution on the initial, unrefined mesh. Each subsequent panel corresponds to a refinement level increased by one, proceeding from top to bottom. Although this is a 3D simulation, only the $x_1x_3$-plane is shown, since the solution is invariant in the $x_2$-direction. The coarsest level has 720 degrees of freedom while the finest has 1.65 million degrees of freedom.}\label{fig:waveguide_solution_Hz}
\end{figure}

Figure \ref{fig:waveguide_solution_Hz} shows a slice of the solution for the $x_3$-component of the magnetic field at different refinement levels for a trial space degree $p=1$. We can see that when the spatial resolution does not respect the Nyquist criterion, the solution does not capture the oscillations of the electromagnetic fields. However, once the critical resolution is reached, the oscillatory pattern is recovered. With about 4 degrees of freedom per wavelength (the third refinement level), the solution shows a decay in the amplitude as we move away from the inlet and a small phase shift towards a smaller frequency. This is an expected behavior \citep[see][]{ihlenburgFiniteElementSolution1995} for under-resolved solutions of wave propagation problems. However, the discontinuous nature of the method localizes the phase shift within the elements, as we see clear jumps in the magnetic field values between adjacent elements. Those observations are in line with other studies using the DPG method for waveguide problems, which noted that the pollution error manifests primarily as an amplitude attenuation \citep{zitelliClassDiscontinuousPetrov2011,HENNEKING202130} and showed that DPG exhibits small phase error in comparison to the standard FEM formulation of time-harmonic wave propagation \citep{demkowiczWavenumberExplicitAnalysis2012}.\par

\begin{figure}
  \centering
  \includegraphics[width=0.8\textwidth]{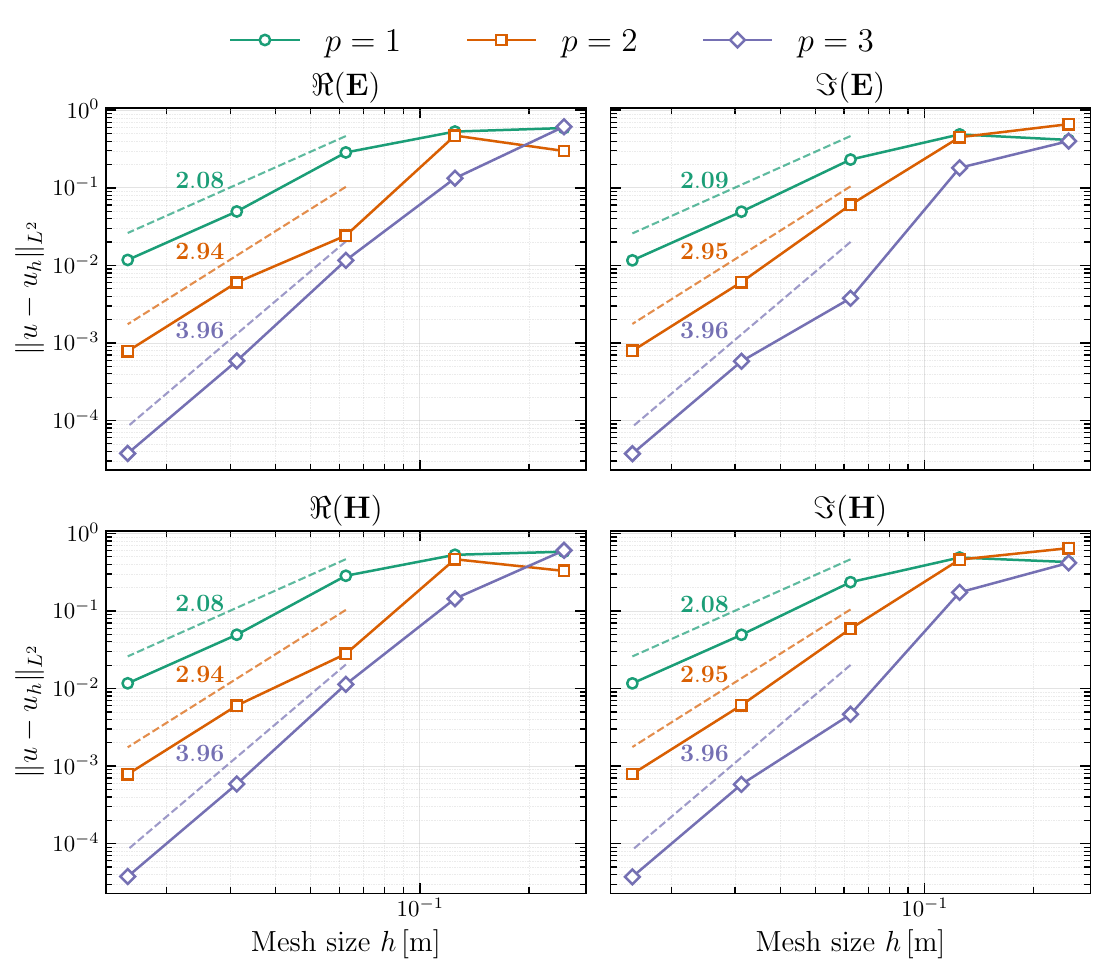}
  \caption{Order of convergence for the real and imaginary part of the electric field $\mathbf{E}$ and the magnetic field $\mathbf{H}$ of interior elements with trial space of order $p\in \{1,2,3\}$ as a function of the mesh size $h$.}\label{fig:waveguide_convergence}
\end{figure}

Figure \ref{fig:waveguide_convergence} shows the convergence of the error for the real and imaginary parts of the electric and magnetic fields as a function of the mesh size $h$ for different trial space polynomial orders. For all field components, we observe the expected optimal convergence rate of $p+1$ for the $L^2$ norm of the error, which is consistent with the theoretical results and confirms the correct implementation of the method. We also note that, for this problem, it is more advantageous to use higher-order polynomials than to refine the mesh. Indeed, the error for the third mesh refinement level with $p=3$ is comparable to the error for the fifth refinement level with $p=1$, while the number of degrees of freedom for the former is about eight times smaller than the latter. This follows from the fact that a 3D hexahedral cell with $p=3$ discontinuous elements $\mathcal{Q}^-_4 \Lambda^3(\square_3)$ has 64 degrees of freedom, but the same cell refined twice (64 cells) with $p=1$ discontinuous elements $\mathcal{Q}^-_2 \Lambda^3(\square_3)$ has 512 degrees of freedom. This is in line with known advantages of higher-order approaches for wave propagation problems in standard FEM (see, for example, \citet{ihlenburgFiniteElementSolution1997}). \par

\subsection{Fichera oven}

\begin{figure}
  \centering
  \includegraphics[width=0.25\textwidth]{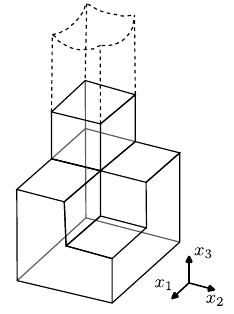}
  \caption{Schematic of the Fichera oven geometry adapted from \citet{carstensenBreakingSpacesForms2016}.} \label{fig:fichera_schematic}
\end{figure}

The second test case is the Fichera oven problem, taken from the DPG literature \citep{carstensenBreakingSpacesForms2016, petridesAdaptiveMultigridSolver2021a}. It assesses the ability of the method to capture singular solutions and to validate the built-in error estimator for adaptive mesh refinement. The geometry of the problem is shown in Figure \ref{fig:fichera_schematic} and consists of a $(0,2)^3$ cube from which a $(0,1)^3$ cube has been removed from one of the corners and appended to the opposite top corner, forming a variation on the 3D L-shape cavity. This design results in a re-entrant corner that introduces a geometric singularity where the electromagnetic fields are expected to exhibit non-smooth behavior, motivating the use of an adaptive mesh refinement strategy.\par

The boundary conditions for this problem are such that all the faces of the cube are perfect electric conductors (PEC) except for the face at the top, which forces an excitation using a Dirichlet boundary condition on the electric field (i.e., $\mathbf{n} \times \mathbf{E} = (\sin{\pi x_2},0,0)$). The problem is expressed in a dimensionless manner. Finally, the properties of the material filling the domain are $\varepsilon_{r,\mathrm{eff}} = 1$ and $\mu_r = 1$, and the frequency of the excitation is $\omega = 5$. \par

\begin{figure}
  \centering
  \includegraphics[width=\textwidth]{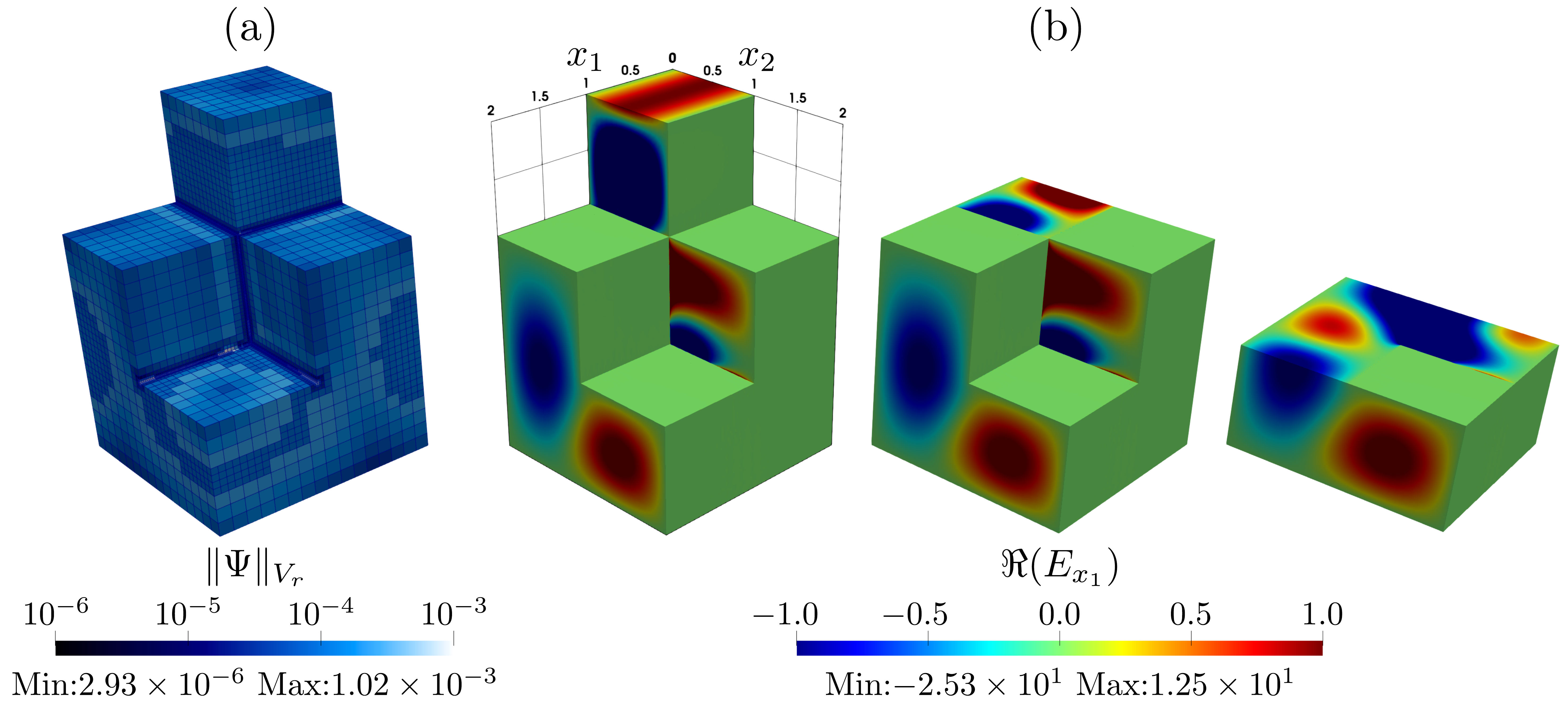}
  \caption{Solution of the Fichera oven test case at the final refinement level with a global error $\|u_h-u\|_E = 0.0232$. In (a) we show the cell residual $\|\Psi\|_{V_r}$, and in (b) we show the $x_1$-component of the electric field. } \label{fig:fichera_solution}
\end{figure}

\begin{figure}
  \centering
  \includegraphics[width=0.5\textwidth]{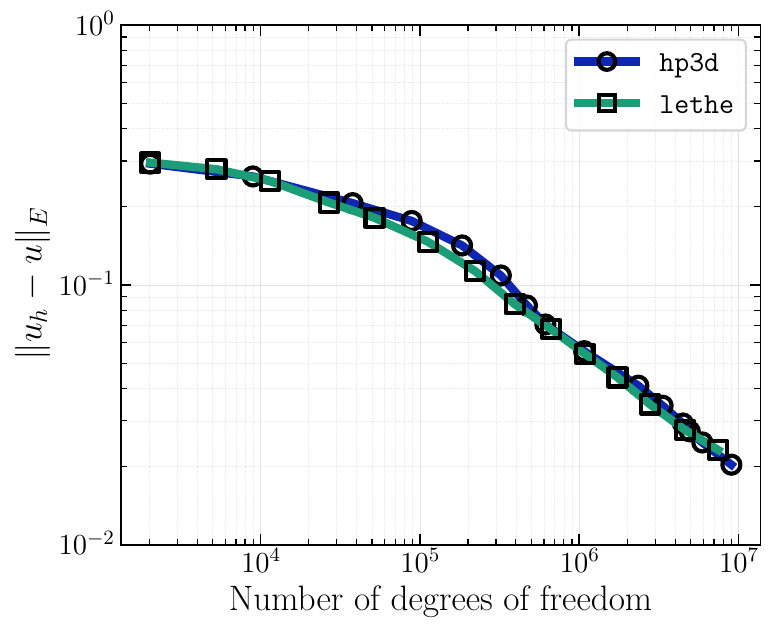}
  \caption{Convergence of the global error in the energy norm $\|u_h-u\|_E$ as a function of the number of degrees of freedom in the linear system for the Fichera oven test case.} \label{fig:fichera_convergence}
\end{figure}

Since there is no known analytical solution for this problem, we reproduce the convergence analysis of the global error $\|u_h-u\|_E$ (Fig. \ref{fig:fichera_convergence}) performed by \citet{petridesAdaptiveMultigridSolver2021a} (Fig. 5.3 (A)). This is obtained by integrating residual norm $\|\Psi\|_{V_r}$ over the domain. In Figure \ref{fig:fichera_solution} (a), we show the cell residual $\|\Psi\|_{V_r}$, and in Figure \ref{fig:fichera_solution} (b), we show the solution of the real part of the $x_1$-component of the electric field at the final refinement level for comparison with the original work.\par

Following \citet{petridesAdaptiveMultigridSolver2021a}, the simulation is initialized with eight cells of trial polynomial degree $p=2$. The mesh is then adaptively $h$-refined based on $\|\Psi\|_{V_r}$ until a number of degrees of freedom comparable to the finest level reported in their study is reached (i.e., approximately $10^7$). The adaptation strategy employed here is based on a Dörfler marking approach, in which cells are selected for refinement based on their contribution to the total estimated error. In the present case, at each adaptation step, cells are ranked according to their residual norms, and the cells with the largest residuals are marked for refinement until their combined residual accounts for 30\% of the total error estimate. Similarly, cells with the smallest residual norms are marked for coarsening until their combined residual accounts for 5\% of the total error estimate.

Analyzing the solution in Figure \ref{fig:fichera_solution} (b), we can see that it is qualitatively similar to the solution of the last refinement level in the original study of \citet{petridesAdaptiveMultigridSolver2021a}, which corresponds to the closest cell residual norm to our final level. Moreover, comparison of the mesh refinement patterns, illustrated by the cell residual in Figure \ref{fig:fichera_solution} (a), shows that, in both studies, the built-in error estimator $\|\Psi\|_{V_r}$ preferentially refines cells in the vicinity of the re-entrant corner and edges, where the solution singularities are located. This behavior provides a strong qualitative validation of the error estimator's ability to correctly identify regions of high error. Nevertheless, a direct comparison with the original study reveals that the final meshes are not identical. This discrepancy is attributable to differences in the mesh refinement strategy, as the refinement procedure employed in \citet{petridesAdaptiveMultigridSolver2021a} is not specified.\par

Similarly, by examining the convergence of the global error $\|u_h-u\|_E$ as a function of the number of degrees of freedom in the linear system presented in Figure \ref{fig:fichera_convergence}, we can see a monotonic decrease in the global error in agreement with the convergence trend of \citet{petridesAdaptiveMultigridSolver2021a}. This further confirms the correct implementation of the method and verifies our $h$-adaptive refinement capabilities.\par

\subsection{Flow around a microwave-heated body}

For the final test cases, we consider the problem of a flow around a microwave-heated body. This case is inspired by the numerical and experimental setup of \citet{pengResonanceDrivenMicrowaveHeating2024}. In their study, the authors investigate the heating of a cylindrical post of different radii placed in the center of a WR430 rectangular waveguide with a cross-section dimension of $x_1=109.2$ by $x_2=\SI{54.6}{\milli\metre}$ and a length of $x_3=\SI{200}{\milli\metre}$. The post is made of alumina ($\varepsilon_{r, \mathrm{eff}} = 9.2 + 0.005i$) and is heated by a microwave excitation at \SI{2.45}{\giga \hertz} with a power of \SI{50}{\watt}. They showed that the system can exhibit resonance-driven heating behavior resulting in a non-monotonic relationship between the post radius and the maximum electric field amplitude, and measured the corresponding temperature increase in the post after 60 seconds of heating.\par

\begin{figure}
  \centering
  \includegraphics[width=0.55\textwidth]{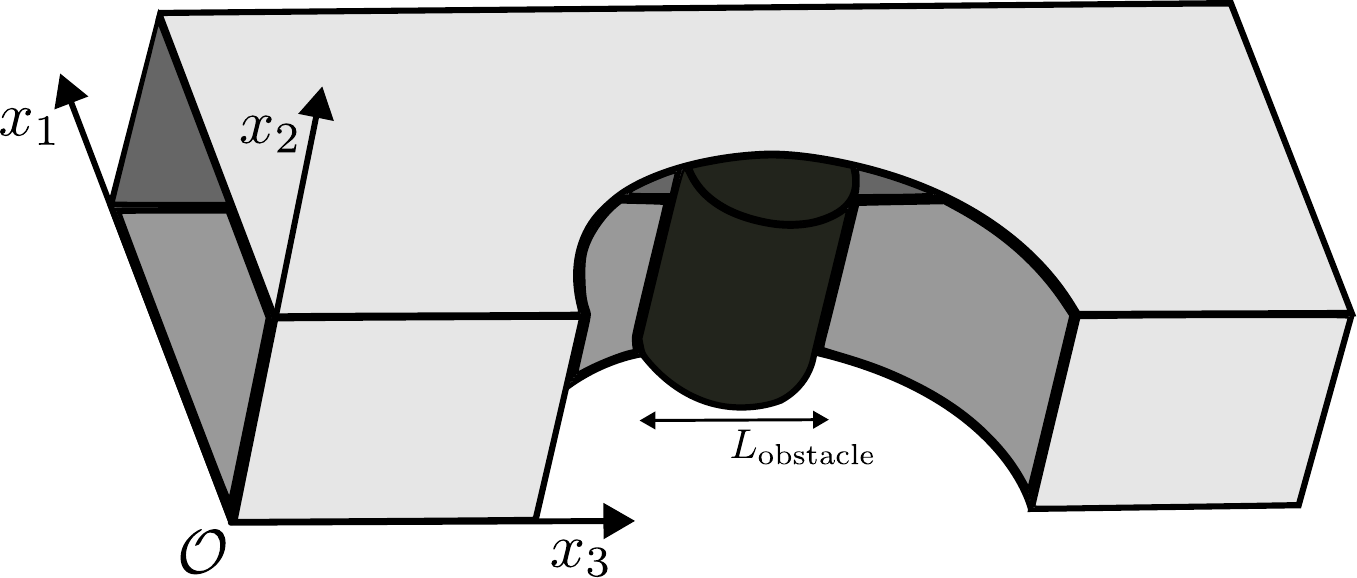}
  \caption{Schematic of the microwave heating setup adapted from \citet{pengResonanceDrivenMicrowaveHeating2024}.} \label{fig:filled_waveguide_schematic}
\end{figure}

Here, we first reproduce the resonance case with a post of radius $\SI{0.024}{\metre}$ to verify the coupled electromagnetic and heat transfer solver against the results of \citet{pengResonanceDrivenMicrowaveHeating2024} without flow. The results of this verification are presented in Section \ref{sec:no_inlet_velocity}. Then, we present new simulation results by extending the case to include a flow around three different obstacle geometries: a cylinder of radius \SI{0.02}{\metre}, a square prism with sides of \SI{0.04}{\metre}, and the same square prism tilted by 45$^{\circ}$ about the $x_2$ axis. We also change the material of the obstacle to silicon carbide ($\varepsilon_{r, \mathrm{eff}} = 9.72 + 2.01i$), a common material used in microwave-heated flow applications due to its high thermal conductivity and dielectric loss. The results of these simulations are presented in Section \ref{sec:with_inlet_velocity}. The physical properties of the materials used across all the simulations are summarized in Table \ref{tab:physical_properties}, and are assumed to be temperature-independent due to a lack of availability for temperature-dependent physical properties.\par

\begin{table}[width=.9\linewidth,cols=4,pos=h]
  \caption{Physical properties for microwave-heated flow test cases \citep{shackelfordCRCMaterialsScience2016, goyalReviewMicrowaveassistedProcess2022,pengResonanceDrivenMicrowaveHeating2024}.}\label{tab:physical_properties}
  \begin{tabular*}{\tblwidth}{@{} LCCC@{} }
    \hline
    \textbf{Parameter}                                                    & \textbf{Air}          & \textbf{Alumina} & \textbf{SiC}   \\
    \hline

    Relative effective permittivity $\varepsilon_{r,\mathrm{eff}}$ \,[--] & $1$                   & $9.2 + i0.005$   & $9.72 + i2.01$ \\
    Relative permeability $\mu_r$ \, [--]                                 & $1$                   & $1$              & $1$            \\

    Thermal conductivity $\kappa$ \,[\si{\watt\per\metre\per\kelvin}]
    & $0.026$               & $26$             & $120.92$       \\
    Specific heat $c_p$ \, [\si{\joule\per\kilogram\per\kelvin}]
    & $1006$                & $1046$           & $27.13$        \\

    Density $\rho$ \,[\si{\kilogram\per\cubic\metre}]
    & $1.225$               & $3750$           & $3210$         \\
    Kinematic viscosity $\nu$ \,[\si{\metre\squared\per\second}]
    & $1.48 \times 10^{-5}$ & –                & –              \\

    \hline
  \end{tabular*}
\end{table}
For the validation case with no inlet velocity, we use the original rectangular domain for which we present a schematic of the geometry in Figure \ref{fig:filled_waveguide_schematic}. For the test cases with an inlet velocity presented in Section \ref{sec:with_inlet_velocity}, the waveguide length is extended to $x_3=\SI{0.4}{\metre}$ and the obstacle is shifted to $x_1=\SI{0.0556}{\metre}$. This $\SI{1}{\milli\meter}$ displacement in the $x_1$-direction is used to accelerate the onset of the vortex shedding behind the obstacle by breaking the symmetry.

The boundary conditions for the different test cases and all physics are presented in Table \ref{tab:BC_heated_flow}. Note that for the case without an inlet velocity, the boundary conditions for the heat transfer problem were not defined in the original study, so here we assume that they are all adiabatic and that the initial temperature of the system is uniform at $T=\SI{293.15}{\kelvin}$. For the test cases with an inlet velocity, we apply a uniform inlet velocity profile with piecewise parabolas on the extremities to respect the no-slip conditions at the walls. Using the volumetric flow rate, the velocity profile amplitude is adjusted so it corresponds to a Reynolds number of Re$=400$ based on the characteristic diameter of the obstacle (i.e., $L_c = L_\text{obstacle} = \SI{0.04}{\metre}$). Additionally, the heat transfer boundary conditions apply a constant temperature $T_{\mathrm{in}}=\SI{293.15}{\kelvin}$ at the inlet, and the domain is initialized at this same temperature. All results later reported are presented as the temperature rise $\Delta T$ relative to this reference.

\begin{table}[width=.9\linewidth,cols=4,pos=h]
  \caption{Boundary conditions for the validation cases without (Case~1) and with (Case~2) inlet velocity. The boundary labels follow the convention used in the waveguide test case presented in Section \ref{sec:waveguide_validation}; headers stand for Electromagnetic (EM), Heat Transfer (HT), and Fluid Dynamics (FD). Note that $Z_\mathrm{s} = (\omega \mu_r)/k_\text{l}$ is the surface impedance of the waveguide port, $k_\text{l}$ is the longitudinal wavenumber of the TE$_{10}$ mode, and $\mathbf{g} = \mathbf{n} \times \mathbf{H}_{\mathrm{TE}_{10}} + Z_\mathrm{s}^{-1} \mathbf{n} \times ( \mathbf{E}_{\mathrm{TE}_{10}} \times \mathbf{n} )$ is the electromagnetic excitation.}\label{tab:BC_heated_flow}
  \begin{tabular*}{\tblwidth}{@{} LCCC@{} }
    \hline
    \textbf{Boundary}     & \textbf{EM}                                                                                                        & \textbf{HT}                     & \textbf{FD}                                                                                         \\
    \hline
    \multicolumn{4}{l}{\textit{Case 1: No inlet velocity}}                                                                                                                                                                                                                             \\
    \hline
    Inlet ($\Gamma_1$)    & $\mathbf{n} \times \mathbf{H} + Z_\mathrm{s}^{-1} \mathbf{n} \times ( \mathbf{E} \times \mathbf{n} ) = \mathbf{g}$ & $\mathbf{n} \cdot \nabla T = 0$ & --                                                                                                  \\
    Outlet ($\Gamma_2$)   & $\mathbf{n} \times \mathbf{H} + Z_\mathrm{s}^{-1} \mathbf{n} \times ( \mathbf{E} \times \mathbf{n} ) = 0$          & $\mathbf{n} \cdot \nabla T = 0$ & --                                                                                                  \\
    Walls ($\Gamma_3$)    & $\mathbf{n} \times \mathbf{E} = 0$                                                                                 & $\mathbf{n} \cdot \nabla T = 0$ & --                                                                                                  \\
    \hline
    \multicolumn{4}{l}{\textit{Case 2: With inlet velocity}}                                                                                                                                                                                                                           \\
    \hline
    Inlet ($\Gamma_1$)    & $\mathbf{n} \times \mathbf{H} + Z_\mathrm{s}^{-1} \mathbf{n} \times ( \mathbf{E} \times \mathbf{n} ) = \mathbf{g}$ & $T = T_{\mathrm{in}}$           & $\mathbf{u} = \mathbf{u}_{\mathrm{in}}$                                                             \\
    Outlet ($\Gamma_2$)   & $\mathbf{n} \times \mathbf{H} + Z_\mathrm{s}^{-1} \mathbf{n} \times ( \mathbf{E} \times \mathbf{n} ) = 0$          & $\mathbf{n} \cdot \nabla T = 0$ & $(-p_{\mathrm{r}}\mathbf{I} + \nu (\nabla \mathbf{u} + \nabla \mathbf{u}^T) ) \cdot \mathbf{n} = 0$ \\
    Walls ($\Gamma_3$)    & $\mathbf{n} \times \mathbf{E} = 0$                                                                                 & $\mathbf{n} \cdot \nabla T = 0$ & $\mathbf{u} = \mathbf{0}$                                                                           \\
    Obstacle ($\Gamma_4$) & --                                                                                                                 & --                              & $\mathbf{u} = \mathbf{0}$                                                                           \\
    \hline
  \end{tabular*}
\end{table}

For all test cases, we perform a mesh convergence analysis for the electromagnetic field amplitude $|\mathbf{E}|$ solution across the heated obstacle center in the $x_1$ and $x_3$ directions at the mid-plane in the $x_2$ direction. For the sake of brevity, results are only shown for the comparison test case with no inlet velocity in Figure \ref{fig:microwave_mesh_convergence}. For the flow cases, mesh convergence for the fluid dynamics is additionally assessed through the lift and drag coefficients:
\begin{equation}
  C_\mathrm{D} = \frac{2 f_{x_3}}{\rho U_\infty^2 L_\text{obstacle}}, \quad C_\mathrm{L} = \frac{2 f_{x_1}}{\rho U_\infty^2 L_\text{obstacle}},
\end{equation}
where $f_{x_3}$ and $f_{x_1}$ are the force per unit of length in the $x_3$ and $x_1$ directions, respectively. Again, for brevity, only the results for the cylinder are presented in Section \ref{sec:with_inlet_velocity} and shown in Figure \ref{fig:cl_profiles}.

\subsubsection{No inlet velocity}\label{sec:no_inlet_velocity}

\begin{figure}
  \centering
  \includegraphics[width=0.75\textwidth]{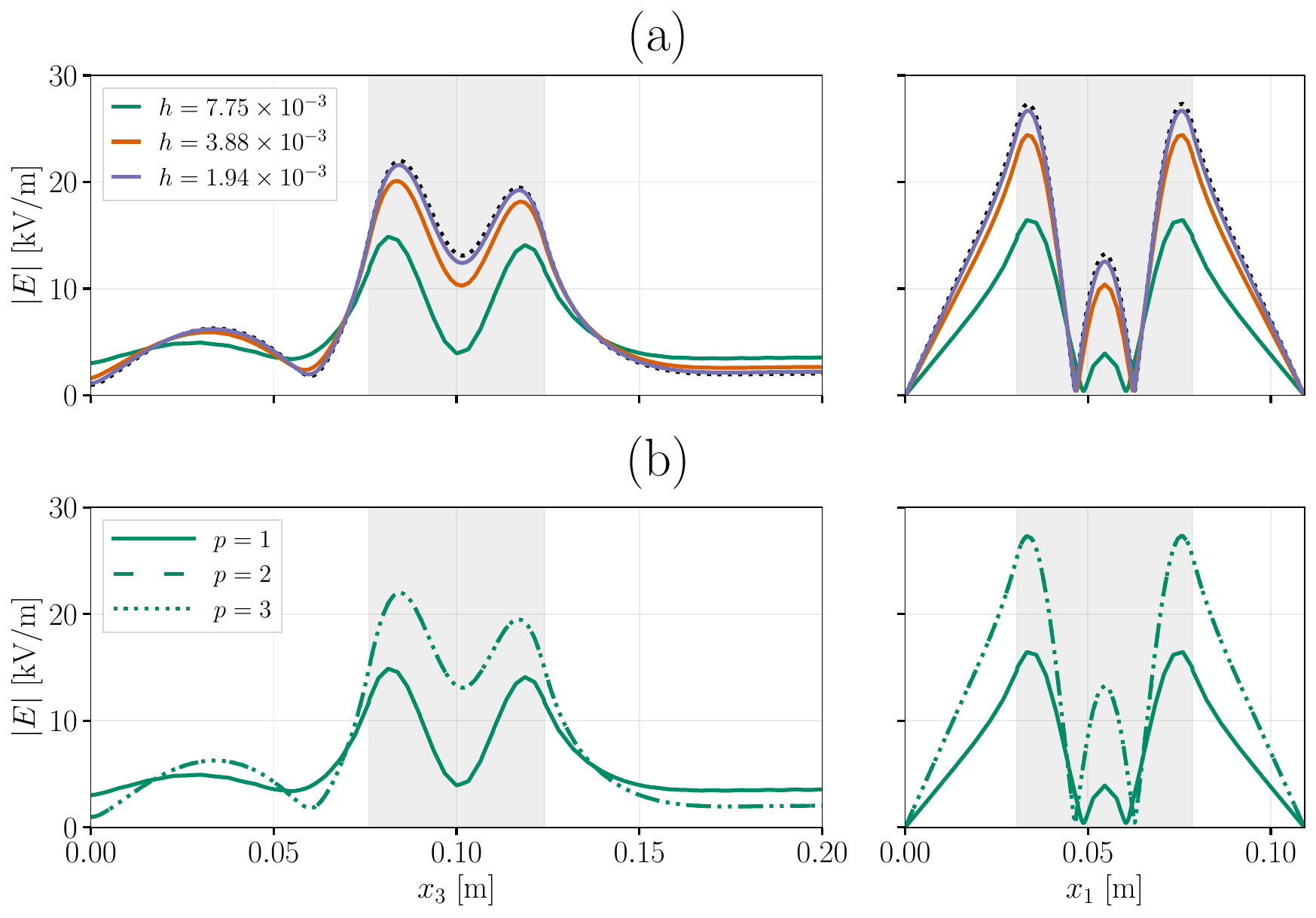}
  \caption{Mesh convergence analysis of the electric field amplitude $|\mathbf{E}|$ solution along the cylinder crossing line at $x_1=0.0546$ and $x_2=0.0273$ m, and the cylinder crossing line at $x_3=0.1$ m and $x_2=0.0273$ m for the test case without inlet velocity. The figures present the convergence for (a) decreasing mean element size $h$ for $p=1$ and (b) increasing polynomial order $p$ for a fixed mean element size of $h=7.75 \times 10^{-3}$ m. The dashed black line in (a) corresponds to the solution obtained with the finest mesh and $p=3$ and is used as a reference. The shaded area represents the cylinder location.} \label{fig:microwave_mesh_convergence}
\end{figure}

\begin{figure}
  \centering
  \includegraphics[width=0.7\textwidth]{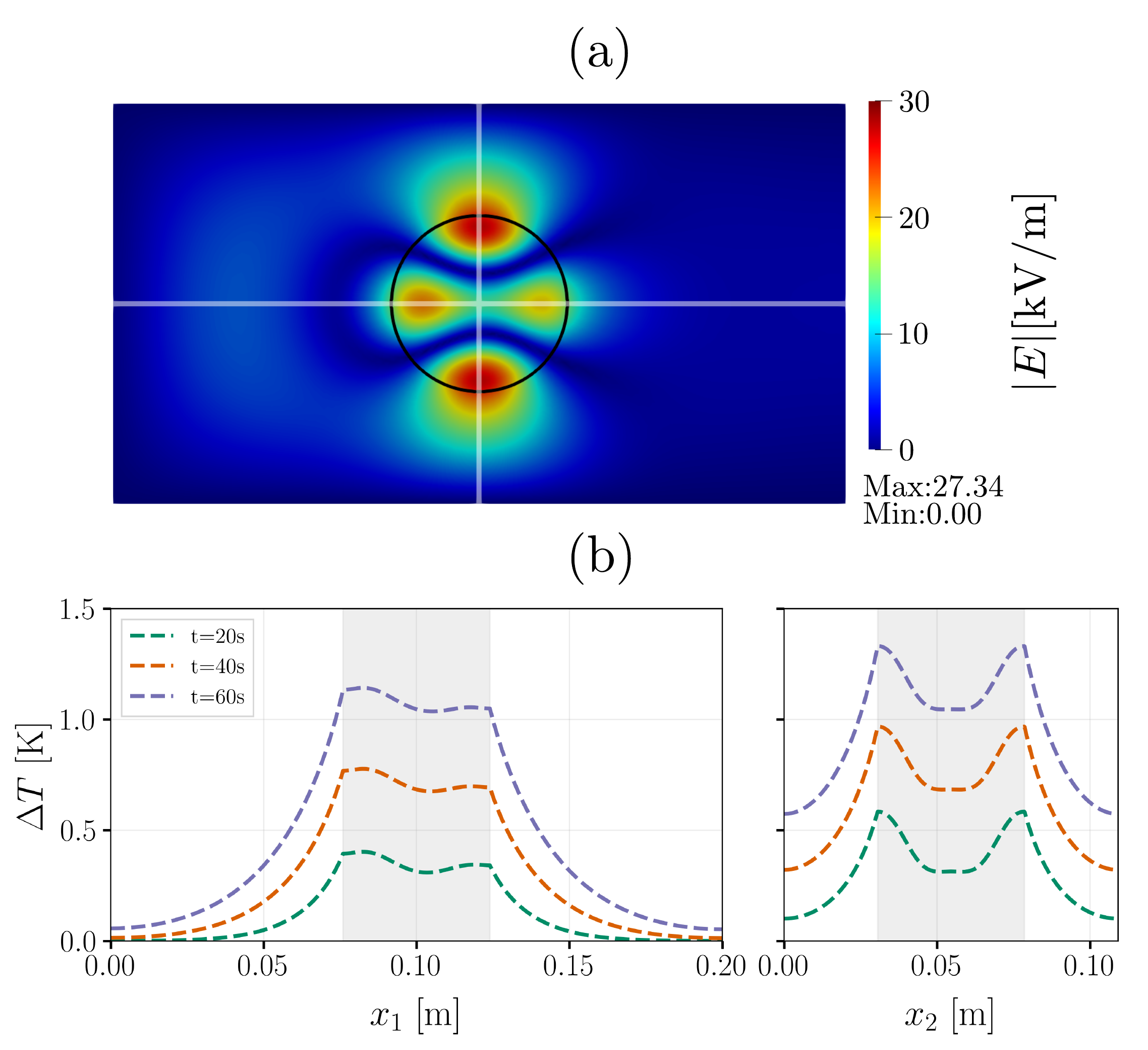}
  \caption{ (a) Electric field amplitude $|\mathbf{E}|$ solution and (b) the temperature profile along the cylinder crossing line at $x_1=0.0546$ and $x_2= \SI{0.0273}{\metre}$, and the cylinder crossing line at $x_3=\SI{0.1}{\metre}$ and $x_2=\SI{0.0273}{\metre}$ at various times.} \label{fig:test_case_comparison}
\end{figure}

From Figure \ref{fig:microwave_mesh_convergence}, we can see that the solution indeed converges as the mesh is refined and as the polynomial order is increased. From it, we can see that the finest mesh for $p=1$ is not fully converged, while the solutions for $p=2$ and $p=3$ are fully converged for a mesh twice as coarse. This reinforces the advantages of using higher-order polynomials for those types of oscillatory problems. The converged electromagnetic field amplitude solution is presented in Figure \ref{fig:test_case_comparison} (a) and shows good agreement with the solution pattern presented by \citet{pengResonanceDrivenMicrowaveHeating2024} (Fig. 12 (c)). However, the amplitude of our solution is slightly lower than the one presented in their study. This is caused by a difference in the definition of the input power for the waveguide port boundary condition, which is a matter of convention. In our case, the input power describes the amplitude of the incident wave that is injected into the domain, while in their case, the input power describes the total power that crosses the inlet plane. This difference in the definition does not affect the validity of the results since the solution pattern is the same and the electric field amplitude can be rescaled by multiplying with the following factor to recover the same amplitude as in their study:
\begin{equation}
  \gamma = \sqrt{P_\mathrm{input}/P_\mathrm{total}} \approx 1.51,
\end{equation}
where $P_\mathrm{total}$ is the total power crossing the inlet plane using the Poynting vector obtained from the solution of $\mathbf{E}$ and $\mathbf{H}$, and $P_\mathrm{input}$ is the input power defined in \eqref{eq:power_scaling}. Note that we do not present the rescaled solution here to be consistent with our definition of the input power and the other results we present later on.

Using an average mesh size of $h=7.75 \times 10^{-3}$ m and a polynomial order of $p=2$ (i.e., the fewest degrees of freedom for a converged solution), we solve the heat transfer problem for 60 seconds of simulation while keeping the electromagnetic solution constant, as no mesh adaptation is performed and all physical properties are assumed constant. The evolution of the change in temperature profile is shown in Figure \ref{fig:test_case_comparison} (b) for different time steps. Note that the white lines in Figure \ref{fig:test_case_comparison} (a) indicate where the temperature and electric field amplitude profiles of (b) and Figure \ref{fig:microwave_mesh_convergence} have been taken, respectively. As expected, the change in temperature profile follows the same pattern as the electromagnetic field amplitude solution, with the temperature increasing more in the regions where the electromagnetic field amplitude is higher. Additionally, without any flow perturbing the system and adiabatic boundary conditions, the spatial shape of the profile is preserved while the post is heated. Finally, after 60 seconds, the average temperature increase in the post is \SI{1.1}{\kelvin}, which is in agreement with the experimental results of \citet{pengResonanceDrivenMicrowaveHeating2024}, where they reported an average temperature increase of approximately \SI{1}{\kelvin} for the same test case. \par

\subsubsection{With inlet velocity}\label{sec:with_inlet_velocity}

From previous studies on the computational efficiency in terms of the memory requirement and time to solution of the matrix-free incompressible Navier-Stokes solver of \texttt{Lethe} \citep{prietosaavedraMatrixfreeStabilizedSolver2025,alphoniusLethe10Opensource2026}, we know that using higher-order polynomials is more efficient than refining the mesh when the solution does not exhibit sharp gradients, as it is expected for the flow around the obstacles at Re$=400$. Therefore, we use a polynomial order of $p=3$ for the fluid dynamics problem and perform the mesh convergence analysis for the lift and drag coefficients using this polynomial order. It follows that the heat transfer problem is also solved with $p=3$ since the Prandtl number of air is close to unity. \par

\begin{figure}
  \centering
  \includegraphics[width=0.6\textwidth]{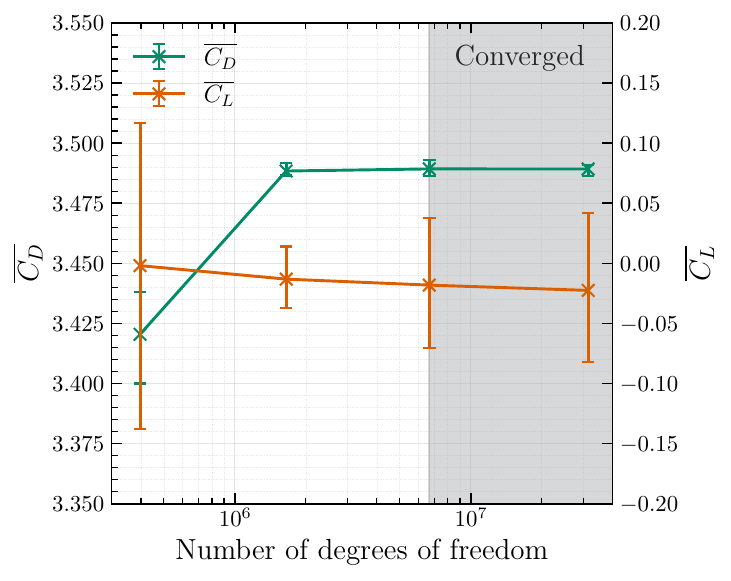}
  \caption{Average in time lift and drag coefficient ($C_\mathrm{l}$ and $C_D$) on the cylinder obstacle as a function of the number of degrees of freedom for the time interval $t\in[6,10]$ \si{\second}. The bar indicates the maximum and minimum values measured in the time interval.} \label{fig:cl_profiles}
\end{figure}

Because the geometry of the problem is rather simple, we use a fixed mesh where additional refinements are concentrated around the obstacle and boundaries to resolve the boundary layers and in the wake of the obstacle. The mesh convergence analysis is then performed by refining the entire domain uniformly. The results of the mesh convergence analysis for the lift and drag coefficients are shown in Figure \ref{fig:cl_profiles}. From it, we conclude that the solution has converged for the third refinement level (i.e., with a mesh with 54368 cells and an average mesh size of $h=3.49 \times 10^{-3}$ m). Now, using the converged solution average element size $h=3.49 \times 10^{-3}$ and the mesh convergence analysis of the electromagnetic field amplitude solution, which showed similar results as the one presented in Figure \ref{fig:microwave_mesh_convergence}, a polynomial order of $p=2$ is sufficient to obtain a converged solution for the electromagnetic field amplitude. Therefore, we use a polynomial order of $p=2$ for the electromagnetic problem and $p=3$ for the fluid dynamics and heat transfer problems for all subsequent simulations. \par

\begin{figure}
  \centering
  \includegraphics[width=\textwidth]{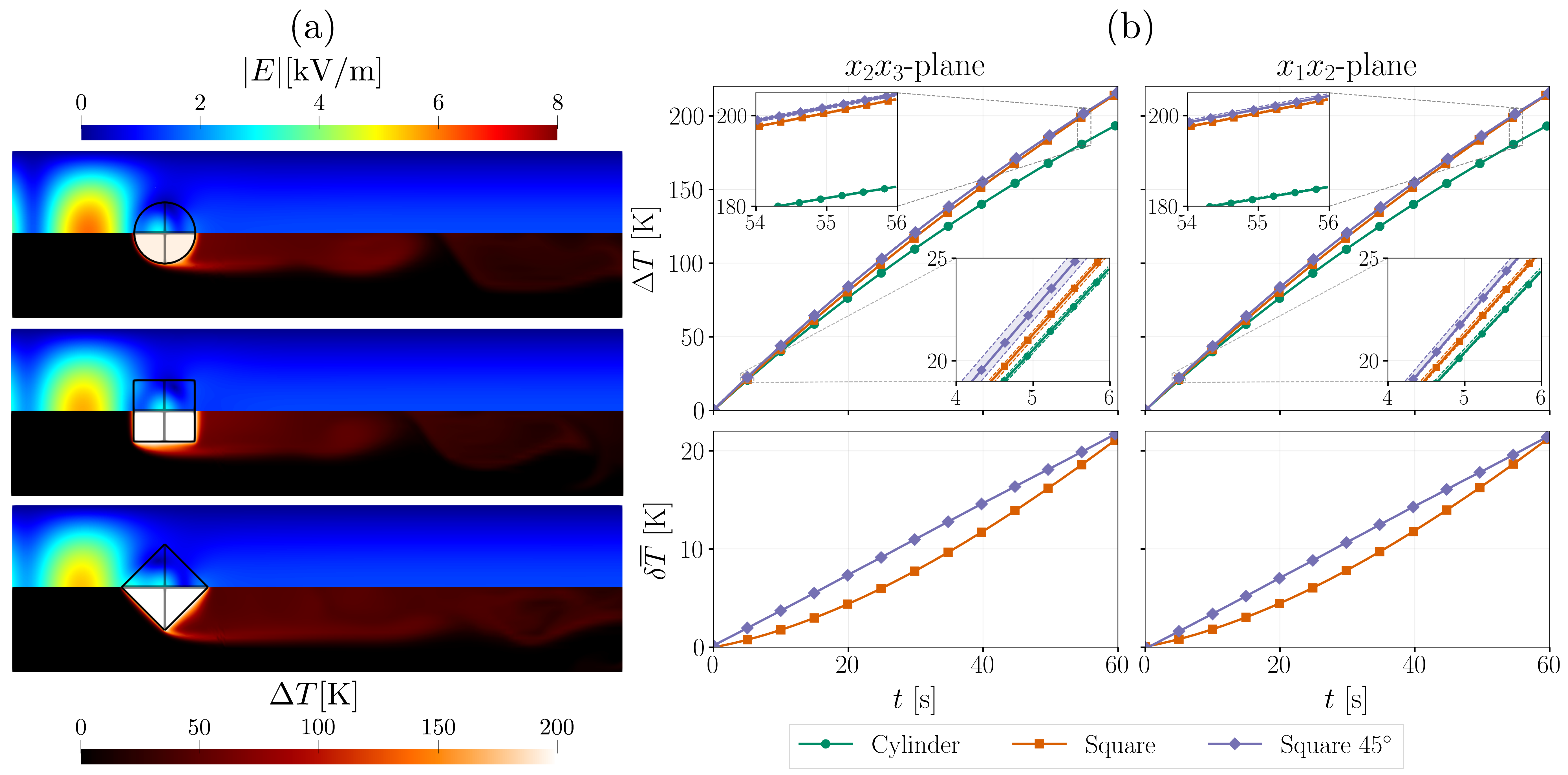}
  \caption{(a) Solution of the temperature and electromagnetic field amplitude after 60 seconds for each different SiC microwave-heated geometry test case. (b) Average change in temperature ($\Delta T$) in the obstacle along the crossing plane at $x_1=\SI{0.0556}{\meter}$ ($x_2x_3$-plane), and the crossing plane at $x_3=\SI{0.1}{\meter}$ ($x_1x_2$-plane) and the difference in average temperature with respect to the cylinder case ($\delta \overline{T}$) along those same planes. The dashed lines bound the region between the $\mathrm{P}_1$ and $\mathrm{P}_{99}$ percentiles. Note that the gray lines in (a) indicate where the profiles of (b) have been taken.} \label{fig:microwave_flow_inlet}
\end{figure}

The simulation results for the flow around the three different obstacles are shown in Figure~\ref{fig:microwave_flow_inlet}. From it, we can see that the electric field amplitude is mostly reflected back toward the incident direction by the obstacle in each case. This situation would be problematic for the microwave-generating source inside a reactor, which would receive the reflected electromagnetic energy and could be damaged without proper shielding. Nonetheless, part of the electromagnetic wave is transmitted through the obstacle and heats it. Note that the penetration depth $D_p$, defined by  \citep{sunReviewMicrowaveMatterInteraction2016}:
\begin{multline}
  D_p = \frac{\sqrt{2}c}{\omega} \bigg( \epsilon_{r,\mathrm{eff},\mathrm{im}}\mu_{r,\mathrm{im}} - \epsilon_{r,\mathrm{eff},\mathrm{re}}\mu_{r,\mathrm{re}} \\
  + \sqrt{(\epsilon_{r,\mathrm{eff},\mathrm{re}}\mu_{r,\mathrm{re}})^2 + (\epsilon_{r,\mathrm{eff},\mathrm{im}}\mu_{r,\mathrm{im}})^2 + (\epsilon_{r,\mathrm{eff},\mathrm{re}}\mu_{r,\mathrm{im}})^2 + (\epsilon_{r,\mathrm{eff},\mathrm{im}}\mu_{r,\mathrm{re}})} \bigg)^{-1/2},
\end{multline}
is approximately \SI{6}{\centi\metre} for the SiC material used for the obstacles. This is similar to the characteristic diameter of the obstacles (i.e., $L_\mathrm{obstacle}= $ \SI{4}{\centi\metre}) and results in an attenuation of the intensity of about 50\%. Therefore, the electromagnetic wave energy is absorbed throughout the obstacle, and only a smaller portion is transmitted through it. Even though the penetration depth and the characteristic diameter are the same for all obstacles, the heating pattern inside each obstacle differs due to their different geometries, which result in different heating rates and temperature distributions, as it can be seen in Figure \ref{fig:microwave_flow_inlet}(b). Indeed, when looking at the zoomed-in views at $t=5$\si{\second} and $t=55$\si{\second}, the tilted square prism shows a more pronounced deviation of its average internal temperature from the dashed lines bounding the $\mathrm{P}_1$--$\mathrm{P}_{99}$ percentile region. This indicates that the tilted square prism has a more heterogeneous temperature distribution than the square prism and the cylinder. This inhomogeneity is, however, on the order of one degree at most and does not change much over time, making it negligible compared to the overall temperature increase later in the simulation. Comparing the two cutting planes also shows that the inhomogeneity is more pronounced in the direction of the flow (which also coincides with the direction of electromagnetic wave propagation) than in the perpendicular direction. This is likely because the flow extracts heat from the obstacle more efficiently at its front than in its wake, resulting in a more pronounced temperature gradient in that direction. In the perpendicular direction, the average temperature is very close to the lower percentile for all geometries, meaning that heat is localized near the obstacle's center and that both sides of the obstacle are rapidly and uniformly cooled by the flow, as expected from the symmetry of the problem. \par

Looking now at the global trend of the average temperature in the obstacle over time ($\Delta T$), we see that the cylinder heats less than the square prism and the tilted square prism, even though the latter two geometries have a volume about 25\% larger than the cylinder. This may seem counterintuitive at first, since more material to heat would suggest a smaller average temperature increase. However, because electromagnetic waves are absorbed volumetrically and the penetration depth is of the same order of magnitude as the obstacles' characteristic diameter, the entire volume of the larger obstacles can absorb electromagnetic energy, resulting in more efficient heating of the larger bodies. After 60 seconds of heating, there is a difference of about \SI{20}{\kelvin} between the average temperature of the cylinder and that of the square obstacles. Additionally, the trend of the average temperature over time shows that thermal equilibrium is not reached in any case for the simulation time prescribed, and that the average temperature is still increasing even after a gain of \SI{200}{\kelvin}. This showcases how efficient microwave heating is compared to conventional heating methods, which would require a much longer time to reach a similar temperature increase. \par

Finally, looking at the difference in average temperature between the square obstacles and the cylinder (i.e., $\delta \overline{T}$), we see that changing only the orientation of an obstacle can yield a different thermal response for the same material and volume. Indeed, the tilted square prism exhibits a heating rate approximately \SI{0.375}{\kelvin\per\second} higher than that of the cylinder, a constant offset that causes the difference between the two temperature profiles to increase linearly over time. The square prism, on the other hand, shows a heating-rate difference relative to the cylinder that evolves in time, starting at approximately \SI{0.2}{\kelvin\per\second} at the beginning of the simulation and increasing to approximately \SI{0.5}{\kelvin\per\second} towards the end. This showcases the nonlinear interplay between the fluid flow and the electromagnetic heating of the obstacle, even under the assumption of constant physical properties. In the current case, it also suggests that at a later time, the square prism would reach a higher temperature than the tilted one, even though the two have similar temperatures at $t=\SI{60}{\second}$.

\section{Conclusions} \label{sec:conclusions}

We have presented an integrated multiphysics software framework that leverages the DPG method to solve the time-harmonic Maxwell equations. The solver can perform adaptive mesh refinement based on the method's built-in error estimator and is coupled to the other physics of the \texttt{Lethe} software, such as heat transfer and incompressible fluid dynamics, to perform multiphysics simulations. All physics support high-order polynomial approximations, the order of which can be chosen independently for each other, allowing for a more efficient use of the degrees of freedom depending on the problem under study. \par

The implementation has been verified against the analytical solution of a rectangular waveguide, since such waveguides are commonly used as the heating source in microwave-assisted reactors. The solver showed the correct convergence order for all the available polynomial degrees implemented (i.e., $p=1,2,3$). It also showed the expected convergence behavior for time-harmonic problems: higher polynomial orders reduce pollution effects, and these effects manifest more as a decrease in solution amplitude than as a complete loss of accuracy \citep{zitelliClassDiscontinuousPetrov2011,HENNEKING202130}. \par

The time-harmonic solver has also been compared against the numerical benchmark of the Fichera oven problem from the DPG literature, showing good agreement with the results of \citet{carstensenBreakingSpacesForms2016,petridesAdaptiveMultigridSolver2021a}. The solver was able to perform mesh adaptation near the domain's singularities and showed a monotonic decrease of the global error in the energy norm as the mesh was refined. \par

Additionally, the solver was verified against the numerical and experimental results of \citet{pengResonanceDrivenMicrowaveHeating2024} for a microwave-heated cylindrical alumina obstacle with a 24 \si{\milli\meter} radius placed in a rectangular WR430 waveguide. The solver reproduced the same electric field amplitude pattern and the same average temperature increase in the obstacle observed after 60 seconds of heating. \par

Finally, we presented new results for three different microwave-heated obstacles simulation with an inlet velocity at a Reynolds number of 400. The results showed how obstacle geometry affects both the internal heating pattern and the surrounding flow responsible for extracting heat from it. They also highlighted the ability of microwave heating to reach high temperatures in a short time, a key advantage of this technology for chemical processing. These final test cases further demonstrated \texttt{Lethe}'s ability to perform multiphysics simulations with different polynomial orders to optimize the use of the degrees of freedom. They also showed the usefulness of the conjugate gradient iterative solver for the time-harmonic Maxwell equations, enabled by the DPG framework, for solving large systems with a limited memory footprint, as required for multiphysics simulations. \par

Future research directions include implementing sum-factorization to reduce the computational cost of the method, an approach already shown to be efficient by \citet{badgerSumFactorizationFast2020}; developing preconditioning strategies to improve the convergence of the iterative solver, building on existing work in this area \citep{petridesAdaptiveMultigridSolver2021a,badgerScalableDPGMultigrid2023}; and extending the solver's coupling capabilities to \texttt{Lethe}'s discrete element method framework for the simulation of microwave-heated multiphase flows \citep{golshanLetheDEMOpensourceParallel2023a,alphoniusLethe10Opensource2026}.

\appendix
\section{Dimensionless conventions}\label{appendix:dimensionless}

The time-harmonic electromagnetic solver is implemented in a dimensionless form using the following scalings:
\begin{align}
   & \mathbf{E} = \frac{1}{E_0}\mathbf{E}_\mathrm{dim}, &                                                     & \mathbf{H} = \frac{Z_0}{E_0} \mathbf{H}_\mathrm{dim},
   &                                                    & \varepsilon_r = \frac{1}{\varepsilon_0}\varepsilon, &                                                                   & \mu_r = \frac{1}{\mu_0}\mu, \nonumber \\
   & \nabla = L \nabla_\mathrm{dim},                    &                                                     & \sigma_r = \frac{1}{\omega \varepsilon_0}\sigma,                  &
   & \omega =  \frac{ L}{c_0}\omega_\mathrm{dim},       &                                                     & \mathbf{J} = \frac{L Z_0}{E_0} \mathbf{J}_\mathrm{dim}. \nonumber
\end{align}
Here, $E_0$ is a reference electric field, $Z_0$, $\varepsilon_0$, $\mu_0$, and $c_0$ are the impedance, permittivity, permeability, and speed of light in vacuum, respectively, and $L$ is the reference length scale of the problem. The subscript "$\mathrm{dim}$" is used to denote which variable has dimension when ambiguous.

\section{Polynomial spaces} \label{appendix:polynomial_spaces}
When implementing the DPG method, care must be taken in the choice of the polynomial spaces for the trial and test functions. For a deeper discussion on this topic, we refer the reader to \citet{carstensenBreakingSpacesForms2016}; here, we will simply restate how the different broken spaces and interface spaces are defined to complement the space definitions presented in section \ref{sec:dpg}. \par
Starting with the standard Sobolev spaces over the domain $\Omega$:
\begin{subequations}
  \begin{align}
     & H^1(\Omega) = \{ u : \Omega \to \mathbb{R}(\text{or }\mathbb{C}) : u, \nabla u \in (L^2(\Omega))^d \},                                                                                   \\
     & H(\text{curl}, \Omega) = \{ \mathbf{E} : \Omega \to \mathbb{R}^d(\text{or }\mathbb{C}^d) : \mathbf{E} \in (L^2(\Omega))^d, \nabla \times \mathbf{E} \in (L^2(\Omega))^d \},              \\
     & H(\text{div},\Omega) =  \{\boldsymbol{\sigma} :\Omega\to\mathbb{R}^d(\text{or }\mathbb{C}^d) :  \boldsymbol{\sigma}\in(L^2(\Omega))^d,  \nabla\cdot\boldsymbol{\sigma}\in L^2(\Omega)\}, \\
     & L^2(\Omega) = \{ q : \Omega \to \mathbb{R}(\text{or }\mathbb{C}) : \|q\|_{L^2(\Omega)} < \infty \},
  \end{align}
\end{subequations}
we define their discretized and broken versions. Limiting ourselves to the discretized 3D case in what follows since our implementation only supports that type of problem, the broken spaces over the domain $\Omega_h$ are defined as:
\begin{subequations}
  \begin{align}
     & H^1(\Omega_h) := \prod_{K \in \Omega_h} H^1(K) = \{ u \in L^2(\Omega) : u|_K \in H^1(K), \forall K \in \Omega_h \},                                                                              \\
     & H(\mathrm{curl}, \Omega_h) := \prod_{K \in \Omega_h} H(\mathrm{curl}, K)=  \{ \mathbf{E} \in (L^2(\Omega))^3 : \mathbf{E}|_K \in H(\mathrm{curl}, K),\ \forall K \in \Omega_h \},                \\
     & H(\mathrm{div}, \Omega_h) := \prod_{K \in \Omega_h} H(\mathrm{div}, K) = \{ \boldsymbol{\sigma} \in (L^2(\Omega))^3 : \boldsymbol{\sigma}|_K \in H(\mathrm{div}, K),\ \forall K \in \Omega_h \},
  \end{align}
\end{subequations}
where $\prod_{K \in \Omega_h}$ is the Cartesian product of the different spaces over all elements $K$ in the mesh $\Omega_h$ and where $|_{K}$ denotes the restriction of the global functions to the element $K$. Note that $L^2(\Omega_h) \equiv L^2(\Omega)$ since $L^2$ functions carry no inter-element regularity requirement and therefore no breaking is needed.\par

Additionally, to define the trace spaces, we first establish the following linear maps and trace operators over the mesh skeleton $\partial \Omega_h$:
\begin{subequations}
  \begin{align}
    \mathrm{tr}_{\mathrm{grad}}        & :
    \begin{cases}
      H^1(\Omega_h) \to \prod_{K \in \Omega_h} H^{1/2}(\partial K), \\
      \mathrm{tr}^K_{\mathrm{grad}}(u) = u|_{\partial K}, \quad u \in H^1(K),
    \end{cases}
    \\[0.6em]
    \mathrm{tr}_{\mathrm{curl},\top}   & :
    \begin{cases}
      H(\mathrm{curl}, \Omega_h) \to \prod_{K \in \Omega_h} H^{-1/2}(\mathrm{curl}, \partial K),                                                                   \\
      \mathrm{tr}^K_{\mathrm{curl},\top}(\mathbf{E}) = ((\mathbf{n} \times \mathbf{E}) \times \mathbf{n})|_{\partial K}, \quad \mathbf{E} \in H(\mathrm{curl}, K), \\
    \end{cases}
    \\[0.6em]
    \mathrm{tr}_{\mathrm{curl},\dashv} & :
    \begin{cases}
      H(\mathrm{curl}, \Omega_h) \to \prod_{K \in \Omega_h} H^{-1/2}(\mathrm{div}, \partial K), \\
      \mathrm{tr}^K_{\mathrm{curl},\dashv}(\mathbf{E}) = (\mathbf{n} \times \mathbf{E})|_{\partial K}, \quad \mathbf{E} \in H(\mathrm{curl}, K),
    \end{cases}
    \\[0.6em]
    \mathrm{tr}_{\mathrm{div}}         & :
    \begin{cases}
      H(\mathrm{div}, \Omega_h) \to \prod_{K \in \Omega_h} H^{-1/2}( \partial K), \\
      \mathrm{tr}^K_{\mathrm{div}}(\boldsymbol{\sigma}) = (\mathbf{n} \cdot \boldsymbol{\sigma})|_{\partial K}, \quad \boldsymbol{\sigma} \in H(\mathrm{div}, K),
    \end{cases}
  \end{align}
\end{subequations}
where $\mathbf{n}$ is the unit outward normal vector to the boundary $\partial K$ of the element $K$, $\mathrm{tr}^K$ denotes the trace operator restricted to the element $K$, and $|_{\partial K}$ denotes the restriction of the global functions to the element boundary. Finally, trace operators are used to define the trace spaces as follows:
\begin{subequations}
  \begin{align}
     & H^{1/2}(\partial \Omega_h) := \mathrm{tr}_{\mathrm{grad}}(H^1(\Omega)),                                   \\
     & H^{-1/2}(\mathrm{curl}, \partial\Omega_h) := \mathrm{tr}_{\mathrm{curl},\top}(H(\mathrm{curl}, \Omega)),  \\
     & H^{-1/2}(\mathrm{div}, \partial\Omega_h) := \mathrm{tr}_{\mathrm{curl},\dashv}(H(\mathrm{curl}, \Omega)), \\
     & H^{-1/2}(\partial \Omega_h) := \mathrm{tr}_{\mathrm{div}}(H(\mathrm{div}, \Omega)).
  \end{align}
\end{subequations}

Finally, the finite element spaces used to discretize the above are denoted, following the finite element exterior calculus convention \citep{arnoldPeriodicTableFinite2014}, by $\mathcal{Q}^-_p\Lambda^k(\square_d)$: the trimmed tensor-product space of degree $p$ for differential $k$-forms on a hexahedral mesh of dimension $d$. On a hexahedral mesh in three dimensions, $\mathcal{Q}^-_p\Lambda^0(\square_3)$ is the space of continuous Lagrange elements, $\mathcal{Q}^-_p\Lambda^1(\square_3)$ that of Nédélec elements of the first kind, conforming in $H(\mathrm{curl})$, $\mathcal{Q}^-_p\Lambda^2(\square_3)$ that of Raviart-Thomas elements, conforming in $H(\mathrm{div})$, and $\mathcal{Q}^-_p\Lambda^3(\square_3)$ that of discontinuous elements, conforming in $L^2$.

\clearpage
\section{Parameters of test cases} \label{appendix:simulation_parameters}
\subsection{Waveguide}
A summary of the solver parameters is presented in Table \ref{tab:params-waveguide-maxwell}.

\begin{table}[width=.9\linewidth,cols=3,pos=h]
  \caption{Parameters for the rectangular waveguide verification simulation.}\label{tab:params-waveguide-maxwell}
  \begin{tabular*}{\tblwidth}{@{} LCC@{} }
    \hline
    \textbf{Parameter}                                 & \textbf{Units} & \textbf{Value}                    \\ \hline

    \multicolumn{3}{l}{\textit{Finite element discretization}}                                              \\ \hline
    Electromagnetic trial degree                       & –              & 1                                 \\
    Electromagnetic test degree                        & –              & 2                                 \\

    \multicolumn{3}{l}{\textit{Physical properties}}                                                        \\ \hline
    Electric conductivity $\sigma_r$                   & –              & 0                                 \\
    Electric permittivity $\varepsilon_r$              & –              & 1                                 \\
    Magnetic permeability $\mu_r$                      & –              & 1                                 \\

    \multicolumn{3}{l}{\textit{Computational domain and mesh}}                                              \\ \hline
    Domain size $(L_x \times L_y \times L_z)$          & \si{\meter}    & $0.25 \times 0.25 \times 1$       \\
    Initial subdivisions $(n_x \times n_y \times n_z)$ & –              & $1 \times 1 \times 4$             \\

    \multicolumn{3}{l}{\textit{Electromagnetic excitation}}                                                 \\ \hline
    Electromagnetic frequency                          & \si{\hertz}    & $2.45 \times 10^9$                \\
    Waveguide mode type                                & –              & TE$_{10}$                         \\

    \multicolumn{3}{l}{\textit{Boundary conditions}}                                                        \\ \hline
    Perfect electric conductor boundary                & –              & ID 0, 1, 2, 3                     \\
    Waveguide port boundary                            & –              & ID 4                              \\
    Impedance boundary                                 & –              & ID 5 : $Z_s^{-1} =  0.968987646 $ \\

    \multicolumn{3}{l}{\textit{Linear solver}}                                                              \\ \hline
    Relative residual tolerance                        & –              & $10^{-8}$                         \\
    Minimum residual                                   & –              & $10^{-12}$                        \\
  \end{tabular*}
\end{table}

\clearpage
\subsection{Fichera oven}
A summary of the solver parameters is presented in Table \ref{tab:params-fichera}.

\begin{table}[width=.9\linewidth,cols=3,pos=h]
  \caption{Parameters for the fichera oven verification simulation.}\label{tab:params-fichera}
  \begin{tabular*}{\tblwidth}{@{} LCC@{} }
    \hline
    \textbf{Parameter}                        & \textbf{Units} & \textbf{Value}                            \\ \hline

    \multicolumn{3}{l}{\textit{Adaptive mesh refinement}}                                                  \\ \hline
    Number of mesh adaptation                 & –              & 13                                       \\
    Error estimator                           & –              & DPG                                       \\
    Refinement fraction                       & –              & 0.30                                      \\
    Coarsening fraction                       & –              & 0.05                                      \\
    Fraction typer                            & –              & fraction                                  \\

    \multicolumn{3}{l}{\textit{Finite element discretization}}                                             \\ \hline
    Electromagnetic trial degree              & –              & 2                                         \\
    Electromagnetic test degree               & –              & 3                                         \\

    \multicolumn{3}{l}{\textit{Physical properties}}                                                       \\ \hline
    Electric conductivity $\sigma$            & –              & 0                                         \\
    Electric permittivity $\varepsilon_r$     & –              & 1                                         \\
    Magnetic permeability $\mu_r$             & –              & 1                                         \\

    \multicolumn{3}{l}{\textit{Computational domain and mesh}}                                             \\ \hline
    Domain size $(L_x \times L_y \times L_z)$ & \si{\meter}    & $2 \times 2 \times 3$                     \\

    \multicolumn{3}{l}{\textit{Electromagnetic excitation}}                                                \\ \hline
    Electromagnetic frequency                 & \si{\hertz}    & $2.38567258 \times 10^8$                  \\

    \multicolumn{3}{l}{\textit{Boundary conditions}}                                                       \\ \hline
    Perfect electric conductor                & –              & ID 0                                      \\
    Electric field boundary condition         & –              & ID 1 : $\mathbf{E} = [\sin(\pi x_2),0,0]$ \\

    \multicolumn{3}{l}{\textit{Linear solver}}                                                             \\ \hline
    Relative residual tolerance               & –              & $10^{-8}$                                 \\
    Minimum residual                          & –              & $10^{-12}$                                \\
  \end{tabular*}
\end{table}
\clearpage
\subsection{Microwave-heated flow}
A summary of the solver parameters for the cylinder microwave-heated body without inlet velocity and with inlet velocity is presented in Table \ref{tab:params-cylinder-resonance} and Table \ref{tab:params-full-coupled-cylinder}, respectively. In the following tables, the acronyms FD, HT, and EM stand for fluid dynamics, heat transfer, and electromagnetism.

\begin{table}[width=.9\linewidth,cols=3,pos=h]
  \caption{Parameters for the transient multiphysics simulation of a microwave-heated cylinder without inlet velocity.}\label{tab:params-cylinder-resonance}
  \begin{tabular*}{\tblwidth}{@{} LCC@{} }
    \hline
    \textbf{Parameter}                        & \textbf{Units} & \textbf{Value}                        \\ \hline
    \multicolumn{3}{l}{\textit{Simulation control}}                                                    \\ \hline
    Time integration scheme                   & –              & BDF1                                  \\
    Final time                                & \si{\second}   & 60                                    \\
    Maximum CFL number                        & –              & 1                                     \\

    \multicolumn{3}{l}{\textit{Mesh}}                                                                  \\ \hline
    Domain size $(L_x \times L_y \times L_z)$ & \si{\meter}    & $0.2 \times 0.1092 \times 0.0546$     \\
    Cylinder diameter                         & \si{\meter}    & 0.048                                 \\

    \multicolumn{3}{l}{\textit{Electromagnetic excitation}}                                            \\ \hline
    Electromagnetic frequency                 & \si{\hertz}    & $2.45 \times 10^{9}$                  \\
    Input power                               & \si{\watt}     & 50                                    \\
    Waveguide mode type                       & –              & TE$_{10}$                             \\

    \multicolumn{3}{l}{\textit{Boundary conditions}}                                                   \\ \hline
    Heat no-flux condition                    & –              & ID 0, 1, 2, 3, 4, 5                   \\
    Perfect electric conductor boundary       & –              & ID 2, 3, 4, 5                         \\
    Waveguide port boundary                   & –              & ID 0                                  \\
    Impedance boundary                        & –              & ID 1 : $Z_s^{-1} = 0.828306014816808$ \\

    \multicolumn{3}{l}{\textit{Finite element discretization}}                                         \\ \hline
    Temperature degree                        & –              & 3                                     \\
    Electromagnetic trial degree              & –              & 2                                     \\
    Electromagnetic test degree               & –              & 3                                     \\

    \multicolumn{3}{l}{\textit{Linear solver}}                                                         \\ \hline
    EM Method                                 & –              & CG                                    \\
    EM relative residual tolerance            & –              & $10^{-4}$                             \\
    EM minimum residual                       & –              & $10^{-8}$                             \\
    HT Method                                 & –              & GMRES                                 \\
    HT Relative residual tolerance EM         & –              & $10^{-3}$                             \\
    HT Minimum residual EM                    & –              & $10^{-6}$                             \\

    \multicolumn{3}{l}{\textit{Nonlinear solver}}                                                      \\ \hline
    HT maximum iterations                     & –              & 10                                    \\
    HT tolerance                              & –              & $10^{-4}$                             \\
  \end{tabular*}
\end{table}

\begin{table}[width=.9\linewidth,cols=3,pos=H]
  \caption{Parameters for the fully coupled transient multiphysics simulation of the microwave-heated cylinder with inlet velocity.}\label{tab:params-full-coupled-cylinder}

  \begin{tabular*}{\tblwidth}{@{} LCC@{} }
    \hline
    \textbf{Parameter}                        & \textbf{Units}    & \textbf{Value}                        \\ \hline

    \multicolumn{3}{l}{\textit{Simulation control}}                                                       \\ \hline
    Time integration scheme                   & –                 & BDF1                                  \\
    Final time                                & \si{\second}      & 60                                    \\
    Maximum CFL number                        & –                 & 1                                     \\

    \multicolumn{3}{l}{\textit{Mesh and refinement}}                                                      \\ \hline
    Domain size $(L_x \times L_y \times L_z)$ & \si{\centi\meter} & $40 \times 10.92 \times 5.46$         \\
    Cylinder diameter                         & \si{\centi\meter} & 4                                     \\

    \multicolumn{3}{l}{\textit{Electromagnetic excitation}}                                               \\ \hline
    Electromagnetic frequency                 & \si{\hertz}       & $2.45 \times 10^{9}$                  \\
    Input power                               & \si{\watt}        & 50                                    \\
    Waveguide mode type                       & –                 & TE$_{10}$                             \\

    \multicolumn{3}{l}{\textit{Boundary conditions}}                                                      \\ \hline
    No-slip boundary                          & –                 & ID 2, 3, 4, 5                         \\
    Outlet boundary                           & –                 & ID 1                                  \\
    Inlet velocity boundary                   & –                 &
    ID~0:\ $[0,0,u_{x_3}]$,                                                                               \\
    &                   &
    $\displaystyle
      u_{x_3}(x_1,x_2)=\mathrm{Re}\,\alpha\,\bigl|f_{x_1}(x_1)f_{x_2}(x_2)\bigr|,$                          \\
    &                   &
    $\displaystyle
      f_\xi(\xi)=
      \begin{cases}
        \xi(\xi-pc_\xi),               & \xi<\dfrac{pc_\xi}{2}, \\[4pt]
        \dfrac{(pc_\xi)^2}{4},         &
        \dfrac{pc_\xi}{2}\le \xi \le c_\xi-\dfrac{pc_\xi}{2},   \\[6pt]
        (\xi-c_\xi)(\xi-c_\xi+pc_\xi), &
        \xi>c_\xi-\dfrac{pc_\xi}{2},
      \end{cases}
    $                                                                                                     \\
    &                   &
    $\displaystyle
      c_{x_1}=10.92,\:
      c_{x_2}=5.46,\:
      \mathrm{Re}=400,\:
      \alpha=1.78212272417098,\:
      p=0.1
    $                                                                                                     \\
    No-flux condition                         & –                 & ID 1, 2, 3, 4, 5                      \\
    Temperature boundary condition            & –                 & ID 0 : $T=0$                          \\
    Perfect electric conductor boundary       & –                 & ID 2, 3, 4, 5                         \\
    Waveguide port boundary                   & –                 & ID 0                                  \\
    Impedance boundary                        & –                 & ID 1 : $Z_s^{-1} = 0.828306014816808$ \\

    \multicolumn{3}{l}{\textit{Finite element discretization}}                                            \\ \hline
    Velocity polynomial degree                & –                 & 3                                     \\
    Pressure polynomial degree                & –                 & 3                                     \\
    Temperature polynomial degree             & –                 & 3                                     \\
    Electromagnetic trial degree              & –                 & 2                                     \\
    Electromagnetic test degree               & –                 & 3                                     \\

    \multicolumn{3}{l}{\textit{Nonlinear solver settings}}                                                \\ \hline
    Fluid dynamics tolerance                  & –                 & $10^{-5}$                             \\
    Heat transfer tolerance                   & –                 & $10^{-4}$                             \\

    \multicolumn{3}{l}{\textit{Linear solver}}                                                            \\ \hline
    FD method                                 & –                 & GMRES                                 \\
    FD relative residual tolerance            & –                 & $10^{-4}$                             \\
    FD preconditioner                         & –                 & GCMG                                  \\
    FD max Krylov vectors                     & –                 & 500                                   \\
    FD coarsening type                        & –                 & PH                                    \\
    FD smoother iterations                    & –                 & 3                                     \\
    FD smoother type                          & –                 & Inverse diagonal                      \\
    FD coarse grid solver                     & –                 & Direct                                \\
    FD eigenvalue estimation                  & –                 & Enabled                               \\
    HT method                                 & –                 & GMRES                                 \\
    HT relative residual tolerance            & –                 & $10^{-4}$                             \\
    EM method                                 & –                 & CG                                    \\
    EM relative residual tolerance            & –                 & $10^{-4}$                             \\
  \end{tabular*}
\end{table}

\printcredits

\section*{Declaration of competing interest}
The authors declare that they have no known competing financial interests or personal relationships that could have appeared to influence the work reported in this paper.

\section*{Acknowledgments}
This research was supported by the Canada Research Chair (CRC-2022-00340), the Natural Sciences and Engineering Research Council of Canada through the RGPIN-2020-04510 Discovery Grant, and computational resources provided by the Digital Research Alliance of Canada. The Trottier Energy Institute also provided financial support through the Hydro-Québec Excellence Scholarship Program, awarded to the first author.

\section*{Data availability}
The data generated and analyzed in this study can be reproduced by running the simulations using the parameters provided in Appendix \ref{appendix:simulation_parameters}. The source code of the \texttt{Lethe} software is openly available at \url{https://github.com/chaos-polymtl/lethe}, with the specific version used to generate the results presented in this article archived at \url{https://doi.org/10.5281/zenodo.22235039}.

\bibliographystyle{cas-model2-names}

\bibliography{cas-refs}

\end{document}